\documentclass[11pt]{article}
\usepackage{amsfonts}
\usepackage{amssymb}
\usepackage{amsmath}
\usepackage{amsxtra}
\usepackage{graphicx}
\usepackage{psfrag}
\usepackage{mathrsfs}
\usepackage{stmaryrd}
\usepackage{subfig}
\usepackage{bm}
\usepackage{tikz} 
\usepackage{empheq}
\usepackage[utf8]{inputenc} 
\usepackage[english]{babel}
\usepackage[square,numbers]{natbib}
\usepackage{pgfplots}
\usepackage{pgfplotstable}
\pgfplotsset{compat=1.3}
\pgfplotsset{every axis legend/.append style={
    at={(1.05,1)},
    anchor=north west,font=\small}}
\usepackage{filecontents}
\usepackage{float}
\usepgfplotslibrary{units}

\usepackage{hyperref}  

\usepackage{xcolor}
\usepackage{comment} 

\title{A numerical investigation of dimensionless numbers characterizing meltpool morphology of the laser powder bed fusion process}
\author{Kunal Bhagat, Shiva Rudraraju \\ \\{\small Department of Mechanical Engineering, University of Wisconsin-Madison, Madison, WI, USA}}
\date{}

\begin{document}
\maketitle 
 
\section*{Abstract}
Microstructure evolution in metal additive manufacturing (AM) is a complex multi-physics and multi-scale problem.  Understanding the impact of AM process conditions on the microstructure evolution and the resulting mechanical properties of the printed part is an active area of research. At the meltpool scale, the thermo-fluidic governing equations have been extensively modeled in the literature to understand the meltpool conditions and the thermal gradients in its vicinity. In many phenomena governed by partial differential equations, dimensional analysis and identification of important dimensionless numbers can provide significant insights into the process dynamics. In this context, a novel strategy using dimensional analysis and the method of linear least squares regression to numerically investigate the thermo-fluidic governing equations of the Laser Powder Bed Fusion AM process is presented in this work. First, the governing equations are solved using the Finite Element Method, and the model predictions are validated by comparing with experimentally estimated cooling rates, and with numerical results from the literature. Then, through dimensional analysis, an important dimensionless quantity - interpreted as a measure of heat absorbed by the powdered material and the meltpool, is identified. This dimensionless measure of heat absorbed, along with classical dimensionless quantities \textcolor{black}{such as P\'eclet, Marangoni, and Stefan numbers}, is used to investigate advective transport in the meltpool for different alloys. \textcolor{black}{Further, the framework is used to study the variations of thermal gradients and the solidification cooling rate.  Important correlations linking meltpool morphology and microstructure evolution related variables with classical dimensionless numbers are the key contribution of this work.}

\section{Introduction}
Additive manufacturing (AM) has proven to be a path-breaking manufacturing paradigm that has the potential to disrupt many of the traditional reductive-type manufacturing processes~\cite {huang2015additive}. A wide variety of AM techniques, suitable for printing metals, glasses, ceramics, and polymers~\citep{gibson2021additive}, are in use today \textcolor{black}{and an optimal AM technique can be selected depending on the material, part complexity and design considerations \cite{konda2017additive}.} \textcolor{black}{Laser Powder Bed Fusion (LPBF) AM process is the focus of this work. }This technique is widely used to print metallic components with intricate geometry to their near-net shape. \textcolor{black}{Components printed using the LPBF process have the potential to exhibit improved material properties as compared to the traditional manufacturing process. In particular, additively manufactured hierarchical stainless steels are not limited by strength-ductility tradeoff unlike traditionally manufactured stainless steel~\cite{wang2018additively}.  Tensile and fatigue properties of additively built Ti-6Al-4V were shown to be superior to mill-annealed Ti-6Al-4V~\citep{qian2016additive}. }However, the properties of the printed components are very sensitive to the choice of the LPBF process parameters and the execution of the printing process.  \textcolor{black}{ Realization of the full potential of AM is not possible unless optimized process parameters can be identified for various alloys used in laser-based additive manufacturing \cite{rittinghaus2022new}. }

The LPBF manufacturing technique is a multi-stage process. Initially, \textcolor{black}{a moving blade of polymer edge (recoater)} spreads a metal powder forming a layer of uniform thickness. A high-energy laser moves over a powder layer bed in a predefined path. This results in a localized melting of the powder metal followed by rapid cooling and fusion of melted powder material on the previously built part. A new layer of the powder is then deposited and the process repeats until the desired part is printed in a layer-by-layer fashion~\citep{gibson2021additive}. This multi-stage additive printing process involves melting and solidification of the material, formation of the localized meltpool, convection cells inside the pool, keyhole formation, improper fusion of the powder, building up of the residual stresses, and sometimes unwanted material deformation, etc.~\citep{king2015laser}. 
Existing literature is focused on understanding the effects of additive process parameters on the properties of experimentally printed components such as the surface roughness of overhang structures~\citep{fox2016effect}, bead geometry and microstructure~\citep{dinovitzer2019effect},  tensile strength~\citep{ning2017additive},  \textcolor{black}{and,  }\textcolor{black}{width and penetration depth of single scan track~\cite{makoana2018characterization}}, etc.  \textcolor{black}{In addition to experimental studies, various LPBF processes, especially meltpool behavior \cite{letenneur2019optimization}, build layers \cite{ mirkoohi2018thermal}, laser heat source \cite{ mirkoohi2019heat} effects have been analytically studied.  \textcolor{black}{ Hybrid modeling that combines analytical models and machine learning-based models is useful in estimating desirable meltpool dimensions and optimized process variables \cite{mondal2020investigation}.}}

On the modeling front for LPBF, literature focused on the modeling of the rich multiphysics aspects of the process has been extensively published.  \textcolor{black}{ Abolhasani \textit{et al.} \cite{abolhasani2019analysis} studied the effect of reinforced materials on the cooling rates and meltpool behavior of AlSI 304 stainless steel using finite element method simulations.} \textcolor{black}{Ansari \textit{et al.} \cite{ansari2019investigation} developed a 3D finite element method based thermal model using a volumetric Gaussian laser heat source to model the thermal profile and meltpool size in selective laser melting process. }\textcolor{black}{The heat diffusion models were reinforced by considering localized dynamic and unsteady fluid flow inside the meltpool.  Dong \textit{et al.} \cite{dong2018effect} considered phase transformation, thermo-physical properties,  heat transfer, and meltpool dynamics in their finite element model to investigate the effect of laser power and hatch spacing on the meltpool. } \textcolor{black}{Similar multiphysics model accounting for heat diffusion and fluid flow was presented by Ansari \textit{et al.} \cite{ansari2021selective} to study the effect of laser power and spot diameter on meltpool temperature in the LPBF process. }\textcolor{black}{Gusarov \textit{et al.}~\citep{gusarov2007heat} focused on heat transfer and radiation physics in their numerical model.} More comprehensive numerical models of the LPBF considers temperature-dependent properties, powdered layer, fluid flow, laser scanning, etc.  \textcolor{black}{Mukherjee \textit{et al.}~\cite{mukherjee2018heat} used comprehensive LPBF numerical models to simulate fluid flow and heat diffusion dynamics for most commonly used alloys.  Khairallah \textit{et al}~\cite{khairallah2014mesoscopic} provided richer insights into LPBF printing of stainless steel using various continuum numerical models. Wang \textit{et al.}~\citep{wang2019powder} coupled finite volume,  discrete element, and volume of fluid methods to rigorously model power spreading,  powder melting, and multi-layer effects during LPBF AM of Ti-6Al-4V alloy.} In trying to capture all the important aspects of the LPBF process in a numerical model, challenges exist in terms of numerous variables, process parameters, and their complex interactions, and \textcolor{black}{these are} outlined in the work of Keshavarzkermani \textit{et al.}~\citep{keshavarzkermani2019investigation} and Fayazfar \textit{et al.}~\citep{fayazfar2018critical}.  

Physical processes with many independent parameters can be analyzed and investigated using dimensional analysis. Traditional areas of physics and engineering, especially fluid mechanics and heat transfer have used dimensional and scaling analysis extensively~\cite {ruzicka2008dimensionless}. Researchers in AM-related problems have recently started incorporating dimensional analysis in their work. Van Elsen \textit{et al.}~\citep{van2008application} provided a comprehensive list of dimensionless quantities that are relevant for the additive and rapid manufacturing process. They justified the usability of the dimensional analysis to investigate complex additive processes like LPBF.  \textcolor{black}{ Dimensionless numbers were shown to assist in choosing previously unknown process parameters for the LPBF process to print Haynes 282, a nickel-based superalloy \cite{islam2022high}. }\textcolor{black}{ Weaver \textit{et. al} \cite{weaver2022laser} demonstrated the application of universal scaling laws to study the effect of process variables such as laser spot radius on the meltpool depth.} Rankouhi \textit{et al.}~\citep{rankouhi2021dimensionless} in their experimental work applied the Buckingham-$\pi$ theorem in conjunction with Pawlowski matrix transformation to present dimensionless quantities that correlate well with the density or porosity of the built component. Their proposed non-dimensional quantity is shown to be applicable across different material properties and machine variables, thereby achieving desirable scaling. Another widely applicable dimensionless quantity called Keyhole number was proposed by Gan \textit{et al}~\citep{gan2021universal}. They made use of dimensionless analysis in conjunction with multiphysics numerical models and high-tech X-ray imagining in their discovery. Keyhole number provides useful insights into the aspect ratio of the meltpool and the transformation of the meltpool from a stable to a chaotic regime. Wang and Liu~\citep{wang2019dimensionless} proposed four sets of dimensionless quantities combining process parameters and material properties. Their experimental analysis shows these numbers can effectively characterize phenomena like the continuity of the track and its size and the part porosity.   \textcolor{black}{ Noh \textit{et al.} \cite{noh2022dimensionless} showed that normalized enthalpy and relative penetration depth in the vertical direction can provide reliable process map for printing single track 3D geometries using selective laser melting process. }

\textcolor{black}{
The  published literature surveyed for this work either uses experimental or numerical methods to propose new dimensionless quantities which are specific to the AM process and are not always related to classical dimensionless numbers used in the fields of fluid mechanics and heat transfer. }Classical dimensionless numbers like the P\'eclet number can provide a good understanding of the impact of process variables on the printed components.  \textcolor{black}{Nusselt, Fourier, and Marangoni number provide a good understanding of the mode of heat transport inside the meltpool for varying laser power and scan speed \cite{ahsan2022global}-\cite{wu2022modeling} }. Cardaropoli\textit{et al.}~\citep{cardaropoli2012dimensional} provided a map for Ti-6Al-4V alloy linking dimensionless quantities with the porosity of built parts. Mukherjee \textit{et al.}~\citep{mukherjee2017dimensionless} used their established numerical models of the LPBF process to simulate the building of the different materials with varied process variables. Using a known set of dimensionless numbers representing heat input, P\'eclet, Marangoni, and Fourier numbers, they made sense of the impact of process parameters on important output variables like temperature field, cooling rates, fusion defects, etc.  

Similar to the meltpool in the LPBF process, the traditional welding process also involves the formation of a weldpool which is the site of various multiphysics interactions and processes. Literature on the use of dimensional analysis to understand the flow patterns in the weldpool offers insights that are relevant to AM. This includes the work by Robert and Debroy~\citep{robert2001geometry} where they highlighted the importance of dimensionless numbers like Prandtl, P\'eclet, and Marangoni in understanding the aspect ratio of the weldpool. Using the numerical models to predict the weldpool shape for a range of materials, they presented the insightful role of these numbers in shaping the weldpool morphology. Similarly, Lu \textit{et al.}~\citep{lu2004sensitivity} also discusses the role of Marangoni convection in affecting the aspect ratio and shape of the weldpool. Their analysis shows that the effect of welding process conditions on the weld part can be understood by looking at the non-dimensional numbers like P\'eclet and Marangoni.  Wei \textit{et al.}~\citep{wei2009origin} showed that the formation of a wavy fusion boundary is linked with the critical values of the Marangoni and Prandtl numbers. Fusion boundaries and shapes have a significant impact on the microstructure of the material.  \textcolor{black}{Asztalos \textit{et al.} \cite{asztalos2022modern} applied dimensional analysis to study the polymer additive manufacturing processes. }
 
As can be seen from the literature review presented, the use of dimensionless numbers to understand the complex interaction of physical processes is gaining attention. However, among the proposed dimensionless quantities, few are universally applicable.  Some of them remain applicable only in the context of a specific study or alloy.  \textcolor{black}{
A universal dimensionless variable or normalized graph can facilitate the comparison of results between different studies using different materials \cite{chia2022process}.} This leaves room for the development of novel approaches to characterize the LPBF process using dimensional analysis.  Our goal in this work is to perform such a dimensional analysis and investigate the relation between meltpool morphology and to a lesser degree, microstructure evolution, and the underlying dimensionless quantities naturally manifested by the thermo-fluidic governing equations of the LPBF process. \textcolor{black}{In this context, a novel numerical strategy is presented here, where the data generated using numerical simulations of the thermo-fluidic model for different alloys and different process parameters was used, along with linear regression analysis, to characterize meltpool morphology in terms of the dimensionless numbers relevant to the meltpool heat and mass transport.}
 
The outline of the paper is as follows: Section~\ref{sec:MathematicalModel} introduces the governing equations of the LPBF process in their dimensional and non-dimensional forms, along with the corresponding numerical formulation suitable to be solved using the Finite Element (FE) method.  Section~\ref{sec:ValidationStudies} covers the validation of our FE-based implementation of the LPBF thermo-fluidic model with experimental results and numerical results from the literature. 
\textcolor{black}{In Section~\ref{sec:NumExperiment},  an empirical analysis based on linear least-squares regression is described to identify an important dimensionless quantity that is interpreted as a measure of heat absorbed by the powdered material and the resulting meltpool.  An important relationship is then identified relating the measure of heat absorbed by the meltpool and classical dimensionless quantities relevant to the thermo-fluidic governing equations of the LPBF.} \textcolor{black}{This is followed by a presentation of simulation results in Section~\ref{sub:Results}, including a discussion on the effects of the dimensionless quantities on the meltpool morphology and the resulting microstructure.} Lastly, concluding remarks are provided in Section~\ref{sec:Conclusions}.

\section{Governing equations of the LPBF process }\label{sec:MathematicalModel}

\subsection{Thermo-fluidic model of the LPBF process}\label{sec:StrongForm} 
\begin{subequations}\label{eq:ThermalFluidModel}
\textcolor{black}{
A discussion of the physical processes underlying LPBF that are relevant to the thermo-fluidic model is outlined in this section.} The schematic in Figure~\ref{fig:SchematicLPBF} shows an outline of the LPBF process. In LPBF, a \textcolor{black}{recoater} spreads a metal powder layer on top of the powder bed or partially built part that is enclosed in an inert environment. A high-intensity laser scans over this powder layer, causing localized melting and fusion of the melted powder on top of the partially built part. At the macro-scale or part-scale, this laser irradiation of the metal powder results in the formation of a {\em meltpool} (\textcolor{black}{also referred to as {\em molten pool} or {\em melting pool} in the literature}) of the liquified metal, that subsequently solidifies. This solidification of the meltpool is driven from the mesoscale, where the liquid melt undergoes a phase transformation to a solid phase, but the solidification is spatially heterogeneous and leads to the formation of dendritic structures and eventually the grain-scale microstructure. The formation of dendrites, their morphology, and the related numerical models have been extensively treated by the authors in a recent publication ~\cite{bhagat2022modeling}.  
\begin{figure}[ht]
  \centering
      \includegraphics[width=0.8\textwidth]{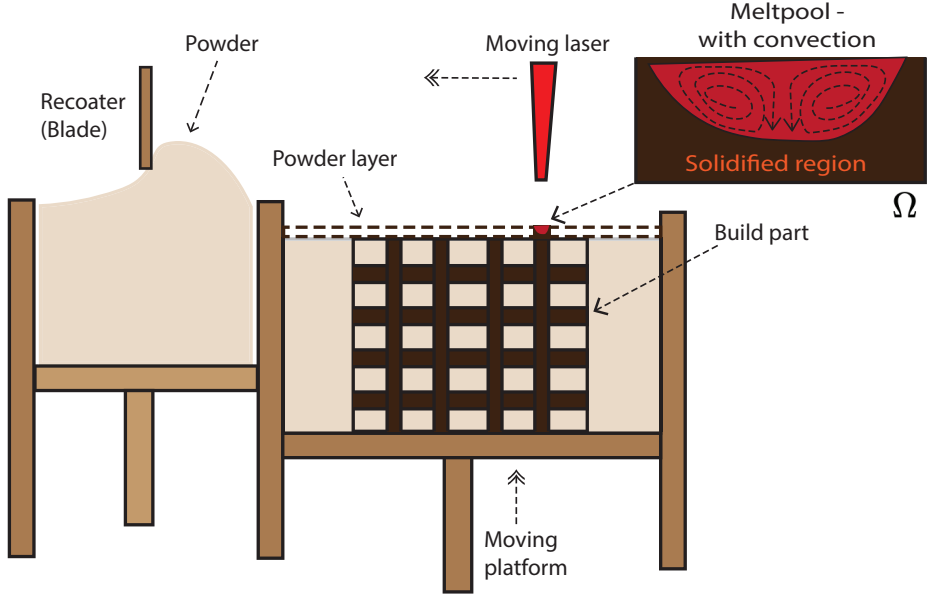}
  \caption{Schematic depicting the part-building process in Laser Powder Bed Fusion (LPBF). Laser irradiation on the powdered material causes localized melting and fusion of the metal powder on top of the partially built part. The localized melting results in a small pool of liquified metal referred to as the {\em meltpool}. Shown in the inset figure are the state of the powder under the laser - with the newly solidified region and a meltpool with convective flow of the liquified metal, and this region comprises the computational domain ($\Omega$) of the numerical model presented in this work.}
  \label{fig:SchematicLPBF}
\end{figure}
\textcolor{black}{
In this work, {\em  part-scale and the thermo-fluidic processes that are relevant in the meltpool and its immediate vicinity are considered.}} The processes \textcolor{black}{modeled}, with varying fidelity, are the movement of the laser-powered heat source, powder melting, convective flow in the meltpool, heat diffusion, and convective and radiation losses. These processes and their thermo-fluidic effects can be captured by coupled partial differential equations solving for the thermal distribution and the velocity distribution in the meltpool~\citep{mukherjee2018heat}. 

\noindent The governing equations of heat transfer are as follows: 
\begin{equation}\label{eq:FluidbasedEnergy}
\rho c \frac{\partial T(\bm{x},t)}{\partial t} + \rho (\bm{\nu}.\bm{\nabla}) T(\bm{x},t) =  \bm{\nabla} \cdot k \big(\bm{\nabla} T(\bm{x},t)\big) +S_{\phi}+S_{p}, \quad  \bm{x} \in \Omega
 \end{equation}
Equation~\ref{eq:FluidbasedEnergy} \textcolor{black}{is solved} for thermal conduction over the domain $\Omega$ (see Figure~\ref{fig:SchematicLPBF}). $T$ and $\bm{\nu}$ are the temperature and velocity, respectively, and are the primal fields of the governing equations. All through this work, vector quantities like velocity and the gradient operator, $\bm{\nabla}$, are shown in bold to distinguish them from other scalar quantities.

$S_p$ is the heat input from the laser and $S_{\phi}$ is the latent heat released by the metal. $\rho$, $c$, $k$ represent density,  specific heat capacity, and the thermal conductivity of the material, respectively, and these can be constant or temperature dependent. Melting of the metal powder consumes latent heat, which is represented by source term, $S_{\phi}=-\rho (\bm {\nu} .\bm{\nabla}) \phi-\rho L (\frac{\partial \phi }{\partial t} )$. Here the liquid fraction, $\phi$, determines the state of the material. $\phi=0$ represents the solidified region,  $\phi=1$ represents the liquid, and $0<\phi <1$ represents the mushy zone region~~\citep{brent1988enthalpy}. The liquid fraction is a function of the temperature of the material and is given by the hyperbolic function: $\phi=\frac{e^{\xi}}{e^{\xi}+e^{-\xi}}$, where $\xi=\frac{\lambda(T-0.5T_{m})}{T_l-T_s}$ and $T_{m}=\frac{1}{2}(T_l+T_s)$.  $\lambda$ is a constant that represents the solid-liquid interface thickness,  $T_s$,  \textcolor{black}{and} $T_l$ are the solidus and the liquidus temperature of the material, respectively. The shape of the laser beam is modeled as an axisymmetric Gaussian profile, and the moving laser power is modeled using a source term, \textcolor{black}{$  S_{p}=\frac{\bm{\alpha} \omega P}{\pi r_{p}^2 l_p} \cdot e^{\big(-\frac{f (x-\bm{\nu}_{p}t)^2}{r_{p}^2}-\frac{f y^2}{r_{p}^2}-\frac{f z^2}{l_p^2}\big)} $}, where $P$ \textcolor{black}{is} the laser power, $\bm{\alpha}$ is the absorptivity,  \textcolor{black}{$f$} is the distribution factor, $r_{p}$ is the laser spot size, $l_p$ is the powder layer thickness, and $\bm{\nu}_{p}$ is the laser scan speed.

\textcolor{black}{Effectively, thermal conduction, powder melting due to the moving laser, and the formation of a liquid meltpool are modeled. }Inside the meltpool, heat diffusion, along with the advection and convection effects of the fluid flow are considered. Convection inside the meltpool is a result of the competition between the surface tension and the buoyancy effects, but the surface tension driven flow \textcolor{black}{dominates} inside the meltpool~~\citep{kumar2009effect}. The governing equations for the fluid flow, accounting for the conservation of mass and momentum, are given by,
\begin{align}
\label{eq:MomentumConservation}
 \rho\frac{\partial \bm{\nu}(\bm{x},t)}{\partial t} +\rho (\bm{\nu}(\bm{x},t).\bm{\nabla}) \bm{\nu}(\bm{x},t) & = -\bm{\nabla} \textcolor{black}{\bm{p}} +\beta \bm{g} (T-T_s) +\bm{ \nabla }. \big(\mu \bm{\nabla} \bm{\nu}(\bm{x},t)\big) -\frac{180\mu}{d_{\phi}^2}\frac{(1-\phi)^2}{\phi^3+\delta}\bm {\nu} (\bm{x},t) \\
 \bm{\nabla}.\bm{\nu}(\bm{x},t) &=0 ,  \quad \bm{x} \in \Omega 
 \end{align}
This equation also accounts for advective and diffusive transport, buoyancy-induced flow, and the pressure drop due to the mushy zone (modeled as a porous zone)~\citep{brent1988enthalpy}. Here, $\beta$ is the expansion coefficient, $T_s$ is the solidus temperature, $\mu$ is the dynamic viscosity, $d_{\phi}$ is the characteristic length scale of the porous mushy region, and $\delta$ is a very small parameter to avoid division by zero when $\phi=0$ (solid region).  \textcolor{black}{Integral form of governing equations introduced in the Section~\ref{sec:WeakFormulation} are given by Equations ~\ref{eq:ND_Energy}-\ref{eq:ND_Momentum}.} As will be seen in the integral form, the surface integrals therein also account for the surface tension-induced flow and losses to the inert surroundings.  
\textcolor{black}{This is captured by the inclusion of the convective and radiation heat losses, Marangoni convection from the top surface of the domain, $\Omega$.  At the bottom surface,  temperature is fixed to a preheating temperature value that is above the ambient temperature. To limit the modeling complexity, in this otherwise highly coupled multiphysics environment, mechanical deformation of the solidified region and keyhole formation is neglected. }

\end{subequations}

\subsection{Non-dimensional formulation of the governing equations}{\label{sec:DimensionlessWeakForm}}
\textcolor{black}{In this section,  a dimensionless framework is constructed that exploits the powerful idea of the Buckingham-$\pi$ theorem.} The LPBF process consists of several process parameters and the thermo-fluidic model \textcolor{black}{that} helps us build an understanding of the complex interaction between several of these process parameters. The dimensionless framework \textcolor{black}{facilitates} combining several of these dimensional parameters into fewer dimensionless independent quantities. These dimensionless parameters then present key details of the complex additive process such as LPBF in fewer numbers of variables. The classical Buckingham-${\pi}$ theorem~\citep{yarin2012pi} provides a mathematical basis for reducing the parameter dimensionality of the equations and helps group the parameters in the governing equations into a fewer number of non-dimensional and distinct quantities. This reduction follows from the application of the Rank-Nullity theorem~\citep{curtis1982dimensional, bluman2013symmetries}.  
\textcolor{black}{Further, as will be discussed in later sections, the Finite Element Method (FEM) is employed to solve the governing equations considered in this work.} FEM is a widely used numerical method that solves partial differential equations posed in their weak formulation (integral formulation). 
\textcolor{black}{Thus, the dimensionless version of the governing equations 
that results from the application of the Buckingham-${\pi}$ theorem is also considered in its weak form and solved using FEM.}

\subsubsection{Weak formulation}{\label{sec:WeakFormulation}}
\textcolor{black}{In this section,  the process of non-dimensionalization of the governing equations given by Equation~\ref{eq:ThermalFluidModel}(a-c) is discussed.  For the process of non-dimensionalization,  the laser scan velocity, $\bm{\nu}_p$ is chosen as the characteristic velocity in the system, and the non-dimensional velocity in the meltpool is then given by $\Tilde{\bm{\nu}}=\frac{\bm{\nu}}{\bm{\nu}_p}$. }The thickness of the powder layer, $l_p$, is chosen as the characteristic length, and this leads to the characteristic time, given by $t_p=\frac{l_p}{\bm{\nu}_p}$.  Now, the non-dimensional time and length are given by $\tilde{t}=\frac{t}{t_p}$ and $\tilde{x}=\frac{x}{l_p}$, respectively. Further, the non-dimensional temperature is chosen to be $\tilde{T}=\frac{T-T_{\infty}}{T_l-T_{\infty}}$, where $T_l$ and $T_\infty$ are the liquidus temperature of the metal and the ambient temperature of the inert surroundings, respectively. The characteristics value of the pressure is chosen to be $\rho \bm{\nu}^2_{p}$. A list of the dimensionless variables used in this model \textcolor{black}{are} summarized in Table~\ref{table:ScalingVariables}. Using these scaled quantities, the dimensional strong (differential) form of the governing equations given by Equation~\ref{eq:ThermalFluidModel}(a-c) are converted into their corresponding dimensionless weak (integral) form. Following the standard variational procedure of deriving the weak formulation of the governing equations from their strong form~\citep{hughes2012finite},  \textcolor{black}{the following weak formulation is obtained:}  \\

\noindent Find the primal fields, $\{ \Tilde{T},\Tilde{\bm{\nu}} \}$, where,
\begin{align*}
\Tilde{T} &\in \mathscr{S}_{T},  \quad \mathscr{S}_{T} = \{ \Tilde{T} \in \text{H}^1(\Omega)~\vert  ~\Tilde{T}  = ~\Tilde{T}' ~\forall ~\bm{x} \in \partial \Omega^{T}_D \}, \\
\Tilde{\bm{\nu}} &\in \mathscr{S}_{\bm{\nu}},  \quad \mathscr{S}_{\bm{\nu}} = \{ \Tilde{\bm{\nu}} \in \text{H}^1(\Omega)~\vert  ~\Tilde{\bm{\nu}} = ~\Tilde{\bm{\nu}}' ~\forall ~\bm{x} \in \partial \Omega^{\bm{\nu}}_D\}
\end{align*}
such that,   
\begin{align*}
\forall ~\omega_T &\in \mathscr{V}_{u},  \quad \mathscr{V}_{T} = \{ \omega_T  \in \text{H}^1(\Omega)~\vert  ~\omega_T   = 0 ~\forall ~\bm{x} \in \partial \Omega^{T}_D \}, \\
\forall ~\bm{\omega_{\nu}}  &\in \mathscr{V}_{\phi},  \quad \mathscr{V}_{\phi} = \{ \bm{\omega_{\nu}} \in \text{H}^1(\Omega)~\vert  ~\bm{\omega_{\nu}}  = 0  ~\forall ~\bm{x} \in \partial \Omega^{\bm{\nu}}_D  \}
\end{align*}
\textcolor{black}{and satisfies,} 

\begin{subequations}
 \begin{align}{\label{eq:ND_Energy}}
  \int_{\Omega} \omega_T \Bigg ( \frac{\partial \Tilde{T}}{\partial \Tilde{t} }
    + (\Tilde{\bm{\nu}}. \Tilde{\bm{\nabla}}) \Tilde{T} \Bigg ) d{{\Omega}}+ 
     \int_{\Omega} \bm{\nabla} \omega_T. 
    \Big[\frac{1}{\bm{Pe}} \Big ] \Tilde{\bm{\nabla}}\Tilde{T}d{{\Omega}} \nonumber 
    +\int_{\Omega} \omega_T \Big [ \frac{\bm {Tc}}{\bm{Ste}}  \Big ] \Big ((\Tilde {\bm{\nu}}. \Tilde{\bm{\nabla}}) \Tilde{\phi} + \frac{ \partial \Tilde {\phi}}{ \partial \Tilde {t}} \Big )d { {\Omega}}\\-\int_{\Omega} \omega_T \frac{\alpha d}{\pi \Tilde{r}^2 \Tilde{l}}  \Big[ \bm{{Q}} \Big]\exp(\tilde{x},\tilde{y},\tilde{z})d{ {\Omega}} 
   + \int_{\partial \Omega^{T}_N} \omega_T \Big (
   \Big [\frac{\bm{Bi}}{\bm{Pe}}  \Big ]
    \Tilde{T}  +
   \Big [ \frac{\bm{t_s}}{\bm{Bo}} \Big ] \Tilde{T}\Big )\bm{n}d{{S}} 
    =0 
    \end{align}
    
\begin{align}{\label{eq:ND_Momentum}}
    \int_{\Omega} \bm{\omega_{\nu}} \Big (\frac{ \partial \Tilde{ \bm{\nu}} }{\partial \Tilde{t}} + \Tilde{\bm{\nu}}\cdot\Tilde{ \bm{\nabla}} \Tilde{\bm{\nu}}\Big )d{{\Omega}} -\int_{\Omega} \bm{\omega_{\nu}} \Big[ \frac{\bm{RaPr}}{\bm{Pe^2}} \Big ] (\Tilde{T}-\Tilde{T_s})d{{\Omega}} - \int_{\Omega} \Tilde{\bm{\nabla}}.\bm{\omega_{\nu}} \Tilde{P}d{{\Omega}} 
    +\int_{\Omega} \bm{\omega_{\nu}} \Big [ \bm{\frac{Pr}{Da Pe}} \Big ]\Tilde{\bm{\nu}} d{{\Omega}} \nonumber \\ +\int_{\Omega}\Big[\frac{\bm{Pr}}{\bm{Pe}}\Big ] (\Tilde{ \bm{\nabla}} \bm{\omega_{\nu}} . \Tilde{\nabla} \Tilde{\bm{\nu}})d{{\Omega}}  +\int_{\partial \Omega^{\bm{\nu}}_N} \bm{\omega_{\nu}}\Big ( \Tilde{P} - \Big [\frac{\bm{Ma Pr}}{\bm{Pe^2}}\Big ]\Tilde{\nabla}  \Tilde{T} \Big )\bm{n}d{{S}} =0
 \end{align} 
\end{subequations}
here, $\bm{n}$ is the unit outward normal vector at the surface boundary. $\partial \Omega^{T}_N$ and $\partial \Omega^{\bm{\nu}}_N$ are the boundaries for the temperature and velocity Neumann boundary conditions, respectively, and $\partial \Omega^{T}_D$ and $\partial \Omega^{\bm{\nu}}_D$ are the boundaries for the temperature and velocity Dirichlet boundary conditions, respectively. $\omega_{T}$ and $\bm{\omega_{\nu}}$ are standard variations from the space $\text{H}^1(\Omega)$ - the Sobolev space of functions that are square-integrable and have a square-integrable derivatives.  In these equations, the relevant dimensionless numbers are grouped inside square brackets. These dimensionless numbers, along with their physical interpretation, are listed in Table~\ref{table:NondimensionalNo}.  The surface boundary condition in Equation~\ref{eq:ND_Energy} represents the nondimensional form of the convective and radiation heat losses to the inert surrounding, and the boundary condition (on the top surface) in Equation~\ref{eq:ND_Momentum} represents the Marangoni flow induced by the surface tension gradient.

\begin{table}[ht]
\centering 
\caption{\textcolor{black}{List of the scaling variables used in the non-dimensionalization of Equations~\ref{eq:ND_Energy}-\ref{eq:ND_Momentum}}} 
\begin{tabular}{c c c} 
\hline 
 Parameter & Expression & Physical interpretation \\ [0.5ex] 
\hline 
$\tilde{l}$ & $\frac{l_p}{l_p}$&  Dimensionless powder layer thickness\\ [1ex]
$\tilde{r}$&	$\frac{r_s}{l}$&	 Dimensionless laser spot radius \\	[1ex]
$\tilde{t}$ & $\frac{t\bm{\nu}}{l_p}$&   Dimensionless time\\ [1ex]
$\tilde{T}$&	$\frac{T-T_\infty}{T_l-T_\infty}$&	 Dimensionless temperature  \\	[1.0ex]
$\tilde{\bm{\nu}}$ & $\frac{\bm{\nu}}{\bm{\nu}_p}$&   Dimensionless velocity \\ [1.0ex]
\textcolor{black}{$\tilde{\bm{p}}$} & $\frac{{\bm{ \textcolor{black}{p}}}}{\rho \bm{\nu}_{p}^{2}}$& Dimensionless pressure \\ [1.0ex]
$\tilde{\bm{\nabla}}$ & $\frac{1}{l_p}\bm{\nabla}$& Dimensionless gradient operator \\ [1.0ex]

\hline 
\end{tabular}
\label{table:ScalingVariables} 
\end{table}


\begin{table}[ht]
\centering 
\caption{\textcolor{black}{Symbols, expressions and their physical interpretation for the dimensionless quantities considered in Equations~\ref{eq:ND_Energy}-\ref{eq:ND_Momentum}}} 
\begin{tabular}{c c c} 
\hline 
 Parameter & Expression & Physical interpretation \\ [0.5ex] 
\hline 
Prandtl $(\bm{Pr})$&	$\frac{\nu }{\alpha}$&	Ratio of momentum to thermal diffusivity \\ [1ex]
Grashof $(\bm{Gr})$&	$\frac{g l^3\beta (T_l-T_{\infty})}{\nu^2}$&	Ratio of buoyancy force to viscous force \\	[1ex]
Darcy $(\bm{Da})$&$\frac{\kappa}{d_{\phi}^2}$& Ratio of permeability to the cross-sectional area\\[1ex]
Marangoni $(\bm{Ma})$&	$\frac{d\gamma}{d T}\frac{l_p\Delta T}{\mu \alpha}$&	Ratio of advection (surface tension) to diffusive transport\\[1ex]
P\'eclet$(\bm{Pe})$&	$\frac{l_p\bm{\nu}_p}{\alpha}$&	Ratio of advection transport to diffusive transport\\	[1.0ex]
Stefan $(\bm{Ste})$& $\frac{c(T_l-T_s)}{L}$&	Ratio of sensible heat to latent heat \\[1.5ex]
Power $(\bm{{Q}})$&	$\frac{P}{\rho c (T_l-T_{\infty})\bm{\nu}_p l_p^2}$&	Dimensionless power with velocity dependence\\[1.5ex]
Radiation measure $(\bm{\frac{t_s}{Bo}})$ &$\frac{\sigma(T_l-T_{\infty})^3}{\rho c \bm{\nu}_p}$& Measure of radiation contribution to the heat transfer\\[1.5ex]
Biot $(\bm{Bi})$&	$\frac{hl_p}{k}$&	Ratio of resistance to diffusion and convection heat transport \\	[1.0ex]
\hline 
\end{tabular}
\label{table:NondimensionalNo} 
\end{table} 
 
\subsection{Computational implementation}
\begin{figure}[ht]
  \centering
      \includegraphics[width=0.55\textwidth]{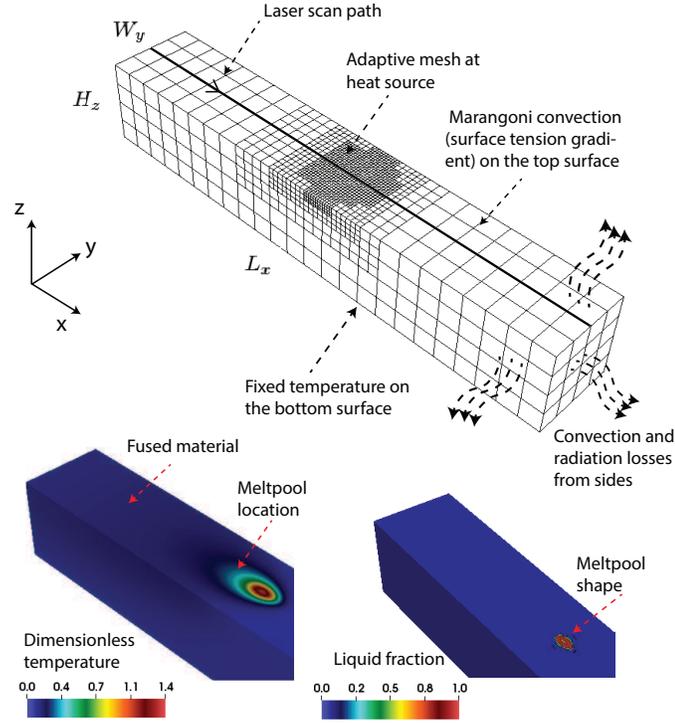}
  \caption{\textcolor{black}{Schematic of the 3D finite element (FE) computational domain indicating the laser scan path and the relevant boundary conditions. Also shown is the underlying adaptive mesh that evolves with the location of the heat source. Representative dimensionless temperature profile and location of the meltpool obtained from the FE simulation of SS316 alloy AM are shown. The numerical parameters and material properties used in this simulation are given in Section~\ref{sec:FEMImplement} and in the Supplementary Information. }}
  \label{fig:FEMSchematic}
\end{figure}

{\label{sec:FEMImplement}}
 \textcolor{black}{As stated earlier,  the above weak formulation of the governing equations is solved using FEM, and as model inputs, realistic process parameters and material properties of common LPBF alloys are chosen,  and these are discussed in \textcolor{black}{Section~\ref{subsec:DataSetGen} and in the Supplementary Information.} 
}FEM implementation is done in an in-house, scalable, finite element code framework written in 
\textcolor{black}{C++ language} with support for adaptive meshing and various implicit and explicit time-stepping schemes, and is built on top of the deal.II open source Finite Element library~\cite{dealII93}.  Standard FEM constructs are adopted, and for all the simulations presented in this work, linear \textcolor{black}{and quadratic Lagrange bases are used for pressure and other field variables such as temperature and velocity, respectively. }The coupled Navier-Stokes equations are solved using Chorin's projection method~\citep{chorin1997numerical}. Following the standard practice in our group to release all research codes as open source~\citep{bhagat2022modeling, GULATI2022100061, wang2016three, jiang2016multiphysics, rudraraju2019computational}, the complete code base for this work is made available to the wider research community as an open-source library~\citep{KunalAdditive2022}. \textcolor{black}{A representative schematic of the computational domain and the relevant boundary conditions are shown in Figure~\ref{fig:FEMSchematic}. The important boundary conditions such as convective and radiations losses and shear stress on the top surface expressed as surface tension gradient with temperature is visible in the surface integral terms in Equation ~\ref{eq:ND_Energy}-\ref{eq:ND_Momentum}.  The initial temperature and temperature at the bottom surface of the material are assumed to be fixed at $353 K$.  The ambient temperature was taken as $301.15 K$.  These temperatures were scaled to a dimensionless form in the computational implementation. The minimum and maximum dimensionless mesh sizes in an adaptive meshing scheme are taken to be $\Delta x = 0.8$ and $\Delta x = 6.0 $ along the x, y, and z directions.  A uniform dimensionless time step size of $\Delta t = 1.0 $ is used for running test cases.  The small factor in Equations~\ref{eq:MomentumConservation},  $\delta = 1.0\times 10^{-5}$ is used in all the simulations. The interface parameter $(\lambda)$ used in our simulation is in the range $\lambda = [0.1,1.0]$}.

\section{Experimental and numerical validation}{\label{sec:ValidationStudies}}
\textcolor{black}{In} \textcolor{black}{this section,  a validation of the numerical formulation and the FEM-based computational framework is presented.  This computational framework solves the thermo-fluidic governing equations of the meltpool described in Section~\ref{sec:MathematicalModel}. }\textcolor{black}{Further,  a comparison is given between the simulation results with experimentally observed cooling rates (made available to us by our experimental collaborators), and with predictions of other numerical models from the literature. }Variables like the cooling rates during the solidification, material temperature, velocity of fluid inside the pool, and meltpool geometry can be used as a yardstick to gauge the capability of our numerical model towards simulating the LPBF process.  \textcolor{black}{For this validation study,  the temperature and velocity distributions, the cooling rates, and the maximum velocity in the meltpool are obtained from our FEM implementation.  } \textcolor{black}{ The cooling rate estimates from our simulations are compared with the cooling rates estimated from experimental data of the LPBF process that were obtained from Bertsch \textit{et al.}~\cite{bertsch2020origin}.  Further,  our simulation results are compared with the corresponding material temperature distribution and meltpool velocity values obtained from numerical modeling data in Shen \textit{et al.}~\citep{shen2020thermo}. } 
Simple thin-walled pseudo-2D plates and 3D cuboidal geometries made of stainless steel (SS316 alloy) using the LPBF process are considered in this validation study.  The printed geometries consists of 13~x~0.2~x~13 $mm^{3}$ thin wall plates and 50~x~10~x~4 $mm^{3}$ cuboids.  \textcolor{black}{The schematic of the printed 2D plates and 3D cuboidal geometries can be found in Bertsch \textit{et al.}~\cite{bertsch2020origin}}.  \textcolor{black}{
These geometries are subsequently referred to as the 2D walls and 3D cuboids.} The powder layer thickness used was 0.02 mm in both cases.  Experimental details, AM technical specifications, and the post-processing methods used to measure cooling rates can be found in the publications of our experimental collaborators, Bertsch \textit{et al.}\cite{bertsch2020origin}-Rankouhi \textit{et al.}\cite{rankouhi2020experimental}.  \textcolor{black}{The experimental cooling rates were estimated by our collaborators, through post-processing of the microstructure morphology, particularly the secondary dendrite arm spacing ($\lambda_2$) at a distance of 6.5 mm and 2 mm from the base for the 2D walls and 3D cuboids, respectively. The dendritic arm spacing in the printed parts was analyzed by our collaborators using a scanning electron microscope (SEM) following standard post-processing techniques. The cooling rates for the alloy SS316 were obtained using the relation $ \lambda_2=25 \dot \epsilon^{-0.28}$ ~\cite{thoma1995directed}, where $\lambda_2$ is measured from SEM images.}

\textcolor{black}{For obtaining the numerical results,  temperature-dependent material properties of the SS316 stainless steel alloy are considered separately for the powdered, fused,  and liquid state of the material.  The temperature and velocity distributions inside the meltpool were obtained from the FE model.} The cooling rates are given by the expression $|\bm{\nabla} T| \bm{\nu}_{p}$, where $|\bm{\nabla} T|_2$ is a measure of the average temperature gradient in the meltpool, and $\bm{\nu}_{p}$ is the laser scanning speed.  For the 2D plates, the cooling rate was measured at a location 6.5 mm from the base, both in the experiments and the FE model.  Similarly, for the 3D cuboids, cooling rate estimates were obtained at a location 2 mm from the base, both in the experiments and the FE model.  As can be seen from Figure~\ref{fig:Experiment_Numerical}, the cooling rates obtained from the FE model are close to the experimentally reported values. 
\textcolor{black}{The cooling rates are used for comparison with experiments in this work,  as they are of immense practical interest due to their influence over the evolution of the microstructure (dendritic growth and grain growth) that then dictates the mechanical properties of the printed part. }

\pgfplotstableread{data/ExpValid.txt}{\ExpValid}

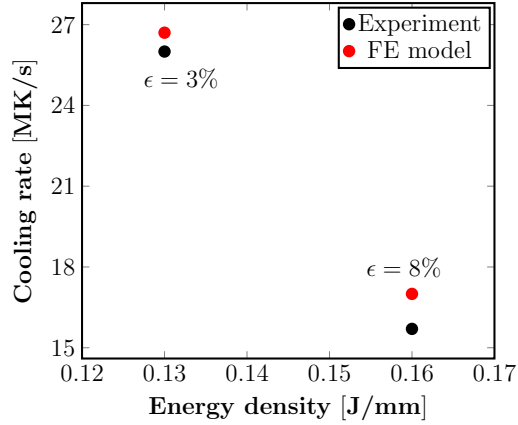
\begin{figure}[h!]
  \centering
    \begin{tikzpicture}[scale=0.65]
        \begin{axis}[ ticklabel style = {thick, font=\Large },very thick,xlabel={\Large {\textbf{Energy density}}},ylabel={\Large \textbf{Cooling rate}},legend style={at={(0.62,0.99)}},ytick={12,15,...,30},x unit=\textbf{J/mm}, y unit=\textbf{MK/s}, mark repeat={1}, scale only axis=true, xmin=0.12,xmax=0.17] 
      \addplot [only marks, black, mark=*, mark size=3] table [x={ED}, y={Exp}] {\ExpValid};
     \addplot [only marks, red, mark=*, mark size=3] table [x={ED}, y={ThmFlud}] {\ExpValid};
    \node[fill =white] at (axis cs:0.132,25.0) {\Large {$\epsilon=3\%$}};
    \node[fill =white] at (axis cs:0.159,18.0) {\Large {$\epsilon=8\%$}};
      \legend{{\Large Experiment},{\Large FE model}}
      \end{axis}
    \end{tikzpicture}  
  \caption{Dependence of cooling rates obtained from experiments and the FE model on the energy density, $\frac{P}{\bm{\nu}_p}$.  The average cooling rate from the FE model was estimated using the relation: $\dot{T}=G\bm{\nu}_{p}=|\bm{\nabla} T| \bm{\nu}_{p}$. Laser power (W) and scan speed (mm/s) combinations used for this study were $(P,\bm{\nu}_p)=(90, 575),  (90, 675)$. } 
\label{fig:Experiment_Numerical}
\end{figure}

Further, the temperature distribution and maximum velocity values in the meltpool obtained from the FE model were compared against the reference test cases given in Shen et al.~\citep{shen2020thermo}. These cases represent the simulation of an LPBF process with a laser power of 100W and 200W used to print AZ91D magnesium alloy parts.  As shown in the Figures~\ref{fig:P200}-\ref{fig:VelP100}, the point temperature as a function of time and the maximum meltpool velocity values obtained from our numerical model closely matches the trend reported in the literature. These comparisons provide a good validation of our FE-based numerical formulation and its computational implementation.  

\pgfplotstableread{data/TempP200.txt}{\TempPowertwo}
\pgfplotstableread{data/TempP100.txt}{\TempPowerone}
\pgfplotstableread{data/VelP100.txt}{\VelPowerone}
\begin{figure}[h!]
  \centering
  \subfloat[\label{fig:P200}]{
    \begin{tikzpicture}[scale=0.5]
 \begin{axis}[ ticklabel style = {thick, font=\Large },very thick,xlabel={\LARGE {\textbf{Time($\times10^{-4}$)}}},ylabel={\LARGE \textbf{Temperature}},legend style={at={(0.35,0.18)}},x unit=\textbf{s}, y unit=\textbf{K}, mark repeat={1}, scale only axis=true] 
    \addplot [only marks, red, mark=*, mark size=3] table [x={Time}, y={Exp}] {\TempPowertwo};
        \addplot [black, mark=none,line width=1.0pt] table [x={Time}, y={Temp}] {\TempPowertwo};   
        \node[fill =white] at (axis cs:1.25,720.0) { \Large {Max $\epsilon=23\%$}};
     \legend{{\Large Shen et. al},{\Large FE Model}}  
      \end{axis}
    \end{tikzpicture}
  }
  \subfloat[\label{fig:P100}]{
    \begin{tikzpicture}[scale=0.5]
     \begin{axis}[ ticklabel style = {thick, font=\Large },very thick,xlabel={\LARGE {\textbf{ Time($\times10^{-4}$)}}},ylabel={\LARGE \textbf{ Temperature}},ytick={0,250,...,1900},legend style={at={(0.35,0.18)}},x unit=\textbf{s}, y unit=\textbf{K}, mark repeat={1}, scale only axis=true] 
    \addplot [only marks, red, mark=*, mark size=3] table [x={Time}, y={Exp}] {\TempPowerone};
        \addplot [black, mark=none,line width=1.0pt] table [x={Time}, y={Temp}] {\TempPowerone};    
        \node[fill =white] at (axis cs:1.35,290.0) { \Large{Max $\epsilon=25\%$}};
     \legend{{\Large Shen et. al},{\Large FE Model}}  
      \end{axis}
    
      \end{tikzpicture}
  }
  \subfloat[\label{fig:VelP100}]{
    \begin{tikzpicture}[scale=0.5]
    \begin{axis}[ ticklabel style = {thick, font=\Large },very thick,xlabel={\LARGE {\textbf{Time($\times10^{-4}$)}}},ylabel={\LARGE \textbf{Velocity}},legend style={at={(0.35,0.18)}},x unit=\textbf{s}, y unit=\textbf{m/s}, mark repeat={1}, scale only axis=true] 
      \addplot [only marks, red, mark=*, mark size=3] table [x={Time}, y={Exp}] {\VelPowerone};
        \addplot [black, mark=none,line width=1.0pt] table [x={Time}, y={Vel}] {\VelPowerone};    
        \node[fill =white] at (axis cs:1.95,0.85) {\Large {Max $\epsilon=9\%$}};
          \legend{{\Large Shen et. al},{\Large FE Model}}  
      \end{axis}
    \end{tikzpicture}
  }
 
\caption{Validation of the FE model results by comparing with corresponding values reported in the literature. ~(\ref{fig:P200}) Variation of point temperature with time for the case P=200W.~(\ref{fig:P100}) Variation of point temperature with time for the case P=100W. ~(\ref{fig:VelP100}) Variation of maximum pool velocity with time for the case P=100W. }
  \label{fig:FEModelvsLiterature} 
\end{figure}
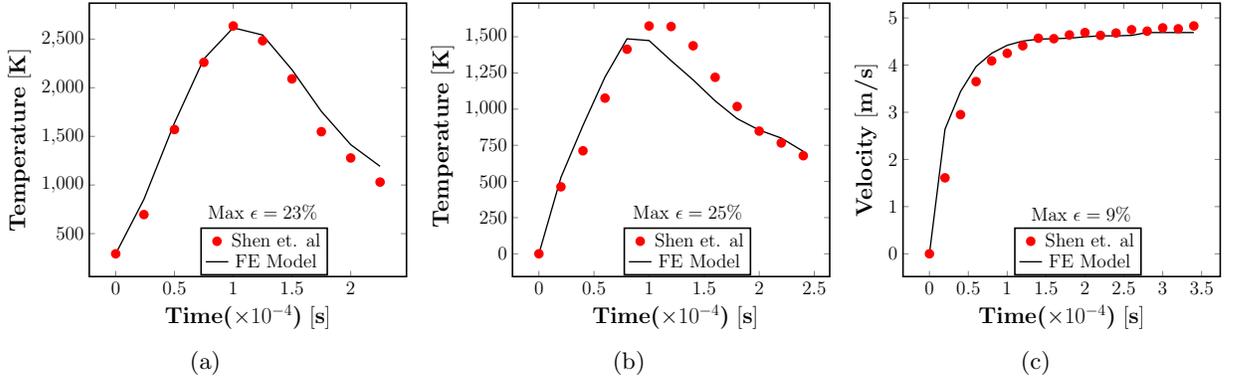

\section{Empirical analysis of the energy absorbed by the meltpool}{\label{sec:NumExperiment}}
\textcolor{black}{In this section,   the rationale behind the construction of a model of the energy absorbed by the meltpool is discussed. }Numerical modeling of the AM process, in general, solves governing equations of heat conduction, fluid flow, mechanical deformation of the solidified region, etc., that is in terms of ``local'' field quantities like temperature, velocity, displacement, etc. However, our goal in this work is to arrive at estimates of the ``global'' response of the system at the meltpool-scale, i.e., measures like the meltpool morphology (size and shape), average temperature distribution, average velocity distribution, etc. {\em The intention is to propose and validate a ``global'' model of the energy absorbed by the meltpool in terms of various process variables, material properties, and output variables, and thus determine the important quantities, from amongst these dependencies, that have a most direct impact on the meltpool evolution}.

Towards proposing a linear model of the heat energy absorbed by the meltpool,  various non-dimensional quantities \textcolor{black}{are chosen} that appear in Equations~\ref{eq:ND_Energy} and \ref{eq:ND_Momentum}. These are the input quantities made up of material properties, processing conditions, and surrounding environment variables. The general procedure used to estimate a linear model is as follows:  \textcolor{black}{$\bm{\hat{U}}$,  the dimensionless measure of the heat absorbed by the meltpool,  and modeled as linearly dependent on the input non-dimensional quantities. }Such a model can be mathematically expressed as  $ \bm{\hat{U}}=\sum_{i=0}^{n}a_i\pi_i$.
where $n$, $\pi_i$, and $a_i$ denote the number of inputs, the dimensionless numbers considered and their corresponding coefficients, respectively. The merits and demerits of choosing various input non-dimensional quantities to characterize the model are evaluated using physics-based and statistical arguments discussed in the subsequent sections. \textcolor{black}{Broadly,  a method of least squares numerical approach is implemented to estimate the coefficients, $a_i$, corresponding to each dimensionless number, $\pi_i$, considered as a potential variable influencing the heat absorbed by the meltpool. } It is the relative magnitude of these coefficients that inform us about the significance or insignificance of each dimensionless quantity towards the model of the heat absorbed.  Further, it is understood from prior knowledge that a system with higher $\bm{\hat{U}}$ can potentially correlate to a rise in some measure of the meltpool temperature.  \textcolor{black}{ The maximum temperature of the material is arguably higher if the heat received $\bm{\hat{U}}$ is higher. }Thus, as a first-order argument, there exists a phenomenological relation $\bm{\hat{U}}\propto \Tilde{T}_{max}$. This reasoning permits us to use $ \Tilde{T}_{max}$ as a measure of the $\bm{\hat{U}}$, and the value of $ \Tilde{T}_{max}$ is obtained by solving thermo-fluidic model described by Equations~\ref{eq:ND_Energy}-\ref{eq:ND_Momentum} on different alloy materials and processing conditions listed in the Tables~\ref{table:AllMaterialProperties}-~\ref{table:AllProcessProperties} provided in the Supplementary Information.  \textcolor{black}{Having obtained $ \Tilde{T}_{max}$,  the coefficients, $a_i$ are estimated, using the method of linear least squares approach and an explicit form $\bm{\hat{U}}$ in terms of various dimensionless numbers, $\pi_i$ is obtained. } Once the linear model of the heat absorbed by the meltpool, $\bm{\hat{U}}$, is determined, it is linked with the several output variables of interest, namely the temperature gradient in the meltpool, the solidification cooling rate $(G\bm{\nu}_p)$,  a measure of the advection heat transport due to the surface tension gradient,  and finally the meltpool morphology (aspect ratio$(\frac{l_{m}}{w_{m}})$ and volume $(l_{m} w_{m} d_{m})$).  Important correlations between the relevant output variables and nondimensional input numbers are discussed at length in Section~\ref{sub:Results}.

\subsection{Process variables, material properties, and output variables }{\label{subsec:DataSetGen}}
\textcolor{black}{In this section, the process variables like laser characteristics, the material properties of the alloy, and the output variables obtained from the thermo-fluidic model are discussed.  The powder layer thickness $l_p=0.02$ mm and laser spot radius $r_p=0.1$ mm are taken for all the cases. }The simulation domain geometry, $L_x\times W_y\times H_z=3\times0.5\times0.5$ $mm^3$, is fixed for all the cases.  
\textcolor{black}{ The movement of the laser is modeled as a single scan on the centerline of the top surface.  The temporal and 3D spatial variations of the temperature and velocity of the material in the meltpool are obtained from the FE model.} Due to the laser heat source, the temperature of the material increases past the liquidus melting temperature and results in the formation of a liquid meltpool.  \textcolor{black}{In the simulations, five commonly used LPBF alloy materials are considered,  namely stainless steel (SS316), a Titanium alloy (Ti-6Al-4V), a Nickel Alloy (Inconel 718), an Aluminium alloy (AlSi10Mg), and a Magnesium alloy (AZ91D)\textcolor{black}{~\cite{khairallah2014mesoscopic,mukherjee2018heat,
shen2020thermo}}.} \textcolor{black}{To limit the complexity of the analysis by making dimensionless quantities independent of temperature,   constant material properties (non-varying with temperature) are chosen. }The numerical values of the input material properties for each of the alloys considered are provided in the Supplementary Information.

The process variables considered are the laser power value $(P)$ and laser scan speed $\bm{\nu}_p$,  \textcolor{black}{a laser distribution factor, $f=2.0$, is fixed for all the cases. }\textcolor{black}{For a given alloy,  twelve $(4\times3)$ combinations of the process variables were chosen to simulate a range of process conditions that are relevant to the LPBF process. }The numerical values of the input process properties for each of the alloys considered are provided in the Supplementary Information, \textcolor{black}{under  Table~\ref{table:AllMaterialProperties} and Table~\ref{table:AllProcessProperties}}.  In the thermo-fluidic model, the heat transfer coefficient $(h)$ and the Stefan-Boltzmann constant $(\sigma)$ are associated with the surrounding inert environment.  \textcolor{black}{ $\sigma=5.67\times 10^{-8}$  $W/(m^2K^4)$ is a known constant.  The effect of varying the heat transfer coefficient is found to have a negligible impact based on our preliminary simulations, so the heat transfer coefficient is taken as $h=10$ $ W/m^{2}K$.} This is due to the relatively minimal convective and radiation losses to the environment, as compared to the conduction of the heat through the base plate. In total, we perform about sixty $(5\times 4\times 3)$ LPBF simulations considering different process variables and material properties.  \textcolor{black}{At a fixed non-dimensional time $\Tilde{t}=100$,  the maximum value of the magnitude of the temperature gradient $G=|\bm{\nabla} T|$ is recorded.  
}The temperature gradient value in the meltpool region is significant but is relatively small outside this region. The temperature gradient is an important variable that controls the microstructure evolution in the additively printed material. The cooling rate, given by $G\bm{\nu}_p$, where $\bm{\nu}_p$ is the speed of the solid-liquid interface is also an important output variable for understanding the microstructure evolution. Further, 
\textcolor{black}{the maximum temperature, $T_{max }$, and maximum velocity, $\bm{\nu}_{max}$, in the meltpool, is tracked along with a measure of the maximum extent of the meltpool length ($l_{m}$), width ($w_{m}$) and depth ($d_{m}$) that characterize the meltpool morphology. }

\subsection{Parametrization in terms of the dimensionless quantities}{\label{subsec:BestLinearModel}}
\textcolor{black}{The use of an empirical approach to estimate~$\bm{\hat{U}}$ is described in this section.  The most appropriate set of dimensionless input parameters that explain variation in the measure of the heat absorbed is chosen.  As stated earlier, $\bm{\Hat{U}}$ is considered proportional to $ \Tilde{T}_{max}$.  Sixty correlations of the form $(\Pi, \Tilde{T}_{max})$  are generated from our simulations, where $\Pi$ represents the set of the dimensionless input quantities considered.  The unknown coefficients, $(a_i)$ are determined through the standard method of linear least-squares regression. } \textcolor{black}{The data for the regression analysis is obtained from the finite element simulations of the LPBF process.} \textcolor{black}{
Here multiple regression attempts were made to arrive at a linear characterization of $\bm{\Hat{U}}$ in terms of the most relevant dimensionless input quantities.  While many combinations of the dimensionless input quantities were considered,  three such attempts as representative of our regression analysis are presented here. }The first two attempts result in correlations that are weak and hence discarded, before converging onto an acceptable correlation between $\bm{\Hat{U}}$ and the relevant dimensionless input quantities in the third attempt.

\paragraph{First attempt of the regression analysis:}The following set of independent variables: $\Pi=\{ \bm{\frac{1}{Pe}},\bm{Q},\bm{\frac{Tc}{Ste}}, \bm{\frac{Bi}{Pe}}, \bm{\frac{t_s}{Bo}}\}$ \textcolor{black}{are considered}. The hypothesized linear relationship is shown below. Here $\epsilon$ is the error - the difference between the expected and true value of $\hat{\bm{U}}$.  
\begin{subequations}
\begin{equation}
\hat{\bm{U}}=a_0+a_1\bm{Q}+a_2\bm{Pe^{-1}}+a_3\bm{\frac{Tc}{Ste}}+a_4\bm{\frac{Bi}{Pe}}+a_4\bm{\frac{ts}{Bo}} + \epsilon
\end{equation}

The values of the coefficients resulting from the least-squares regression are given in Table~\ref{table:OLS_Model0}. The condition number for this analysis is $8.11\times10^7$, which is very high. This indicates that there exists strong collinearity in the assumed input set $\Pi$. The collinearity can be understood in terms of the primary variable as follows: \textcolor{black}{
On close inspection of the dimensionless expressions for $\bm{Q}=\frac{P}{\rho c (T_l-T_{\infty})\bm{\nu}_p l_p^2}$, $\bm{\frac{1}{Pe}}=\frac{\alpha}{\bm{\nu}_p l_p}$, $\bm{\frac{Bi}{Pe}}= \frac{hl\alpha}{\bm{\nu}_p l_p}$, and $\bm{\frac{t_s}{Bo}}= \frac{\sigma (T_l-T_{\infty})^3 } {\rho c \bm{\nu}_p}  $,  it is observed that laser scan velocity is accounted for in all the four variables, and this can potentially reduce the linear independence of these physical quantities. }Further, the role of inert environment variables is limited in our analysis. Considering the regression coefficients,  $\bm{Bi/Pe}$ and $\bm{t_s/Bo}$ \textcolor{black}{are dropped} from our next regression attempt.

\begin{table}[h]
\centering 
\caption{\textcolor{black}{First attempt of the regression analysis to estimate the coefficients, $a_i$, using the linear least-squares approach. Asterisk($^*$) indicates the statistical significance of the coefficient using a t-test with a 95\% confidence interval. Other statistics : $R^2=0.65$, Adjusted $R^2=0.61$, F-statistic=$20.12$, P(F)=$0.0$. Condition number=$8.11\times10^7$.}} 
\begin{tabular}{c c c c c c c c} 
\hline 
Parameter & Intercept & $\mathbf{Q}$ & $\mathbf{Pe^{-1}}$&$\mathbf{\frac{Tc}{Ste}}$&$\mathbf{\frac{Bi}{Pe}}$ &$\mathbf{\frac{t_s}{Bo}}$\\ [0.5ex] 
\hline 
$a_{i}$ & $1.45^{*}$ & $0.0053^{*}$ &$-0.1719^{*}$&$-0.7076^{*}$&$300.3$& $-18200^{*}$\\[0.5ex] 
\hline 
\end{tabular}
\label{table:OLS_Model0} 
\end{table}
\begin{table}[h] 
\centering 
\caption{\textcolor{black}{Second attempt of the regression analysis to estimate the coefficients, $a_i$, using the linear least-squares approach.  Asterisk($^*$) indicates the statistical significance of the coefficient using a t-test with a 95\% confidence interval. Other statistics: $R^2=0.763$, Adjusted $R^2=0.750$, F-statistic=$55.95$, P(F)=$0.0$. Condition number=$1.12\times10^3$. }} 
\begin{tabular}{c c c c c} 
\hline 
Parameter & Intercept& $\bm{E}$ & $\bm{Pe}$&$\bm{\frac{Tc}{Ste}}$\\ [0.5ex] 
\hline 
$a_{i}$ & $0.6938^{*}$ & $0.0087^{*}$ & $-0.1677^{*}$ & $0.1521$ \\[0.5ex] 
\hline 
\end{tabular}
\label{table:OLS_Model1} 
\end{table}
\begin{table}[h]
\centering 
\caption{\textcolor{black}{Third attempt of the regression analysis to estimate the coefficients, $a_i$, using the linear least-squares approach. Asterisk($^*$) indicates the statistical significance of the coefficient using a t-test with a 95\% confidence interval. Other statistics : $R^2=0.746$, Adjusted $R^2=0.737$, F-statistic=$83.61$, P(F)=$0.0$. Condition number=$422$. 
}} 
\begin{tabular}{c c c c} 
\hline 
Parameter & Intercept& $\bm{E}$ & $\bm{Pe}$\\ [0.5ex] 
\hline 
$a_{i}$ & $0.8146^{*}$ & $0.0082^{*}$ & $-0.1654^{*}$ \\[0.5ex] 
\hline 
\end{tabular}
\label{table:OLS_Model2} 
\end{table}

\paragraph{Second attempt of the regression analysis:} \textcolor{black}{
In the second iteration, the chosen independent set is $\Pi=\{\bm{E},\bm{Pe}, \bm{\frac{Tc}{Ste}}\}$. }The hypothesized linear relationship is given by the following relation. 
\begin{equation}
\hat{\bm{U}}=a_0+a_1\bm{E}+a_2\bm{Pe}+a_3\bm{\frac{Tc}{Ste}} + \epsilon
\end{equation}
Further, considering the details of the regression analysis,  \textcolor{black}{
the expression for the variable dimensionless power, $\bm{Q}$, that had a velocity dependence is modified to a new variable, $\bm{E}=\frac{P}{k(T_l-T_\infty)}=\bm{PeQ}$. }This new variable is a modified dimensionless power and is independent of the laser scan velocity. Now, the laser scan speed parameter is only represented in the P\'eclet  number, $\bm{Pe}$. The fourth term $\bm{\frac{Tc}{Ste}}$ is purely dependent on the material properties.  \textcolor{black}{Note} the least-squares solution yields the coefficients given in Table~\ref{table:OLS_Model1}.

The least-square solution obtained from this model is an improvement over the previous model. This can be realized from the improvement in the accuracy of the fit given by the variable $R^2$. The adjusted $R^2$ improves from $0.65$ to $0.75$ with the less number of variables in the set $\Pi$. The condition number is still high but better than the previous model. Thus, the non-significant variable \textcolor{black}{is dropped} from the set, $\Pi$, in our third attempt. 

\paragraph{Third attempt of the regression analysis:}\textcolor{black}{$\Pi=\{\bm{E},\bm{Pe}\}$ and the hypothesized linear relationship is given by the following relation. }
\begin{equation}\label{eg:bestmodel}  
   \hat{\bm{U}}=a_0+a_1\bm{E}+a_2\bm{Pe} + \epsilon  
\end{equation} 
The least-square solution \textcolor{black}{is summarized in} Table~\ref{table:OLS_Model2}. The condition number is greatly improved. The probability that all $a_i=0$ at the same time is negligible as seen from the probability of F-statistic. All the coefficients are statistically significant i.e the hypothesis that individual $a_i=0$ can be safely discarded.

Now the linear model given by Equation~\ref{eg:bestmodel} \textcolor{black}{is interpreted} in light of the physics of the additive process. The heat received by the material is defined by a quantity $\hat{\bm{U}}$. The higher the heat received greater the temperature reached by the system. This $\hat{\bm{U}}$ is related to the various material property and process parameters.  \textcolor{black}{Using the non-dimensional analysis,  several parameters were combined into a bunch of non-dimensional numbers. }These numbers are associated with and are responsible for several physical phenomena.  \textcolor{black}{
Using the linear model and the available data,  the variation of $\hat{\bm{U}}$ was explained using a linear combination of constant, $\bm{E}$ and $\bm{Pe}$. }A more complete dependence can potentially be highly non-linear, but this also can be analytically intractable. 
\end{subequations}
 
\section{Results}{\label{sub:Results}}
\textcolor{black}{In the previous section,  a linear least-squares regression approach was used to arrive at a relation for the dimensionless heat energy absorbed, $\hat{\bm{U}}$. }The regression analysis resulted in a relation for $\hat{\bm{U}}$ in terms of the dimensionless power, $\bm{E}$, and the P\'eclet number, $(\bm{Pe})$.  \textcolor{black}{ In the following sections (Sections~\ref{subsec:AdvectionMeasure}-\ref{subsec:VolumeMeasure}), relation for $\hat{\bm{U}}$ is used to investigate the advective transport occurring inside the meltpool for different alloys} and then \textcolor{black}{to} characterize their meltpool morphology (aspect ratio and volume) using the Marangoni number and the Stefan number.  In Section~\ref{subsec:MicroEvolu},  $\hat{\bm{U}}$ \textcolor{black}{is used} to characterize microstructure evolution using the temperature gradients and the cooling rates in the solidified region. 


\subsection{Influence of P\'eclet number on advection transport in the meltpool}{\label{subsec:AdvectionMeasure}}
\textcolor{black}{In this section, the extent of advection transport observed in the meltpool in different alloys during the LBPF process is discussed. 
The} goal is to analyze the macroscopic geometric features of the meltpool, such as its aspect ratio and volume. In trying to explain the variation of these macroscopic features, a measure of advection in the meltpool using the relevant dimensionless quantities \textcolor{black}{are critically investigated}. Figure~\ref{fig:AdvectMeasure} describes the variation of  the P\'eclet number, $\bm{Pe^{*}}=\bm{Pe\nu_{max}}$, with the non-dimensional quantity $\bm{Ma\hat{U}}$. It is to be noted that the P\'eclet number with an asterisk, $\bm{Pe^{*}}=\bm{Pe\nu_{max}}=\frac{l_p\bm{\nu}_p}{\alpha}\frac{\bm{\nu}_{max}}{\bm{\nu}_p}$, is defined using the maximum velocity in the meltpool, and is a measure of the advective transport relative to the diffusion transport in the meltpool.  A larger value of $\bm{Pe^{*}}$ denotes a larger circulation of heat inside the meltpool due to the fluid motion. The dimensionless quantity $\bm{Ma^*}=\bm{Ma\hat{U}}$ is a measure of heat transport caused by the fluid flow induced due to the surface tension gradient.  As seen from the Figure~\ref{fig:AdvectMeasure}, for the alloy shown, $\bm{Pe^*}$ increase with the $\bm{Ma^*}$.  This correlation implies the overall movement of fluid in the meltpool is greater if the advection transport due to surface tension force is greater.  Each point in these plots represents a single simulation result for the relevant quantities plotted and is obtained from the FEM framework. Another key information conveyed in Figure~\ref{fig:AdvectMeasure} is that for some alloys like AlSi10Mg, advection due to surface tension forces is minimal, as can be seen from the numerical values of the total advection ($\mathbf{Pe\bm{\nu}_{max}}$) shown on the Y-axis.  On the other hand, alloys like Ti6Al4V show a higher value of advection transport due to surface tension forces. These observations are important correlations between advection measure $\bm{Pe^{*}}= \bm{Pe\nu_{max}}$,  Marangoni number, P\'eclet number, and the dimensionless power $(\bm{Ma\hat{U}}=a_0\bm{Ma} +a_1\bm{MaE}+a_2\bm{MaPe})$, and will be used below to make connections to the meltpool morphology.

\pgfplotstableread{data/DimlessQ.txt}{\DimlessQ}
\pgfplotstableread{data/AdvectMeasureSS.txt}{\AdvectMSS}
\pgfplotstableread{data/AdvectMeasureTi.txt}{\AdvectMTi}
\pgfplotstableread{data/AdvectMeasureAl.txt}{\AdvectMAl}

\begin{figure}[ht]
  \centering
\subfloat[\label{AdvectAl}]{
   \begin{tikzpicture}[scale=0.5]
      \begin{loglogaxis}[ ticklabel style = {thick, font=\Large },very thick,xlabel={\Large { $\bm{Ma\hat{U}}$ }},ylabel={\Large \textbf{$\bm{Pe}\bm{\nu}_{max}$}},legend style={at={(0.68,0.11)}},x unit=, y unit=, mark repeat={1},scale only axis=true] 
     \addplot [only marks, cyan, mark=*, mark size=3] table [x={AL_MaUhat}, y={AL_PeU0}] {\AdvectMAl};   
   \addplot [thick, cyan] table[y={create col/linear regression={y=AL_PeU0}}]{\AdvectMAl};  
      \legend{{\Large AlSi10Mg}}
      \end{loglogaxis}
    \end{tikzpicture}
  }
     \subfloat[\label{AdvectSS}]{
   \begin{tikzpicture}[scale=0.5]
      \begin{loglogaxis}[ ticklabel style = {thick, font=\Large },very thick,xlabel={\Large { $\bm{Ma\hat{U}}$ }},ylabel={\Large \textbf{$\bm{Pe}\bm{\nu}_{max}$}},legend style={at={(0.77,0.11)}},x unit=, y unit=, mark repeat={1},scale only axis=true] 
       \addplot [only marks, red, mark=*, mark size=3] table [x={SS_MaUhat}, y={SS_PeU0}] {\AdvectMSS};
        \addplot [thick, red] table [x={Logx}, y={FitLogy}] {\AdvectMSS};
             \legend{{\Large SS316}}
      \end{loglogaxis}
    \end{tikzpicture}
  }
    \subfloat[\label{AdvectTi}]{
   \begin{tikzpicture}[scale=0.5]
      \begin{loglogaxis}[ ticklabel style = {thick, font=\Large },very thick,xlabel={\Large { $\bm{Ma\hat{U}}$ }},ylabel={\Large \textbf{$\bm{Pe}\bm{\nu}_{max}$}},legend style={at={(0.71,0.11)}},x unit=, y unit=, mark repeat={1},scale only axis=true] 
       \addplot [only marks, black, mark=*, mark size=3] table [x={TI_MaUhat}, y={TI_PeU0}] {\AdvectMTi};
          \addplot [thick, black] table[y={create col/linear regression={y=TI_PeU0}}]{\AdvectMTi};  
             \legend{{\Large Ti6Al4V}}
      \end{loglogaxis}
    \end{tikzpicture}
  }
  \caption{ Measure of total advection measured as $\mathbf{Pe\bm{\nu}_{max}}$ vs surface tension based advection $\mathbf{Ma\hat{U}}=a_0\bm{Ma} +a_1\bm{MaE}+a_2\bm{MaPe}$ on a log-log scale for (\ref{AdvectAl})AlSi10Mg,  (\ref{AdvectSS})SS316, (\ref{AdvectTi})Ti6Al4V alloys.  \textcolor{black}{Corresponding plots comparing IN718 and AZ91D alloys, and a comparison of all the five alloys considered in this work can be found in Figure~\ref{fig:AdvectMeasureSupplementary} and Figure~\ref{fig:AdvectMeasureFull} of the Supplementary Information, respectively.} The advection measure corresponds to the degree of fluid flow inside the meltpool.  Each point in these plots represent a single simulation result for the relevant quantities plotted, and is obtained from the FEM framework.}
  \label{fig:AdvectMeasure} 
\end{figure}
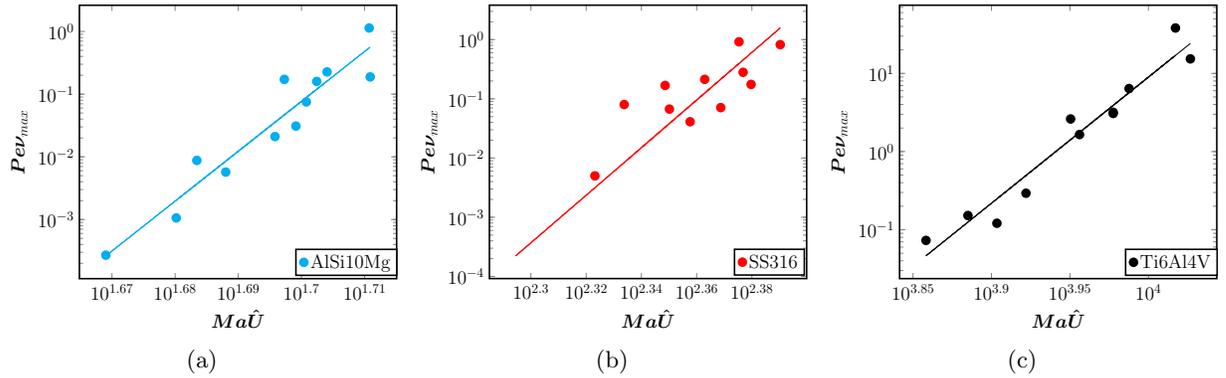

\subsection{Influence of Marangoni number on the meltpool aspect ratio}{\label{subsec:AspectRatioMeasure}}
\textcolor{black}{In this section,  the geometric characteristics of the meltpool, particularly the aspect ratio is discussed.  The meltpool aspect ratio is defined as the ratio of its maximum length to its maximum width $(\frac{l_{m}}{w_{m}})$.  To understand aspect ratio in terms of the input process parameters, dimensionless quantity $\bm{Ma\hat{U}}$ is useful.  For characterizing the meltpool shapes,  correlations were found between $\bm{Ma\hat{U}}= a_0\bm{Ma}+a_1\bm{MaE}+a_2\bm{MaPe}$ and the aspect ratio of the meltpool for different alloys, and different process parameters. } The rise in $\hat{\bm{U}}$ is an indication of enough heat received by the material to cause melting during the additive process. The aspect ratio of the meltpool is related to the fluid flow induced inside the meltpool.  As shown in the previous section, in Figure~\ref{fig:AdvectMeasure}, AlSi10Mg alloy has a low advection measure, $\bm{Pe\nu_{max}}$, and Ti6Al4V alloy has a high advection measure.  This information is relevant here to understand the meltpool shapes of these materials.  Figure~\ref{fig:AspectPrMaUhatTi} shows the variation of the aspect ratio $(\frac{l_{m}}{w_{m}})$ plotted as a function of non-dimensional quantity $\bm{Ma\hat{U}}$ for Ti6Al4V alloy material.  The aspect ratio increases with  $\bm{Ma\hat{U}}$.  \textcolor{black}{ The aspect ratio with $\bm{Ma\hat{U}}$ if visualized in a combined plot for all the three alloys considered in this work,  it is instructional to see the separations of materials into three clusters - each for one alloy, as seen in the Figure~\ref{fig:AspectPrMaUhatCombined}.  
From this clustering, it can be seen that Ti6Al4V almost always produced an elongated or elliptical-shaped meltpool whose aspect ratio is far from one. The alloy AlSi10Mg produces a meltpool that is relatively less elongated and has an aspect ratio closer to one. } The advection in the fluid flow causes the meltpool to expand along the direction of the higher temperature gradient. From this discussion,  \textcolor{black}{insightful observations, relating the aspect ratio of the meltpool with $\bm{Ma\hat{U}}$, can be made that help us characterize the meltpool shapes potentially produced by different alloys. }



\pgfplotstableread{data/AspectRatioSS.txt}{\AspectUhatSS}
\pgfplotstableread{data/AspectRatioTi.txt}{\AspectUhatTi}
\pgfplotstableread{data/AspectRatioAl.txt}{\AspectUhatAl}

\begin{figure}[ht]
  \centering    
    \subfloat[\label{fig:AspectPrMaUhatTi}]{
    \begin{tikzpicture}[scale=0.55]
      \begin{loglogaxis}[ ticklabel style = {thick, font=\Large },very thick,xlabel={\LARGE { $\bm{{Ma\hat{U}}}$ }},ylabel={\huge\bm{$\frac{l_{m}}{w_{m}}$}},legend style={at={(0.71,0.1)}},x unit=, y unit=, mark repeat={1},scale only axis=true] 
         \addplot [only marks, black, mark=*, mark size=3] table [x={TI_MaU}, y={TI_AR}] {\AspectUhatTi};
       \legend{{\Large Ti6Al4V}}
      \end{loglogaxis}
    \end{tikzpicture}
    }   
  \subfloat[\label{fig:AspectPrMaUhatCombined}]{
    \begin{tikzpicture}[scale=0.55]
      \begin{loglogaxis}[ ticklabel style = {thick, font=\Large },very thick,xlabel={\LARGE { $\bm{{Ma\hat{U}}}$ }},ylabel={\huge\bm{$\frac{l_{m}}{w_{m}}$}},legend style={at={(0.675,0.25)}},x unit=, y unit=, mark repeat={1},scale only axis=true] 
         \addplot [only marks, black, mark=*, mark size=3] table [x={TI_MaU}, y={TI_AR}] {\AspectUhatTi};
       \addplot [only marks, red, mark=*, mark size=3] table [x={SS_MaU}, y={SS_AR}] {\AspectUhatSS}; 
    \addplot [only marks, cyan, mark=*, mark size=3] table [x={AL_MaU}, y={AL_AR}] {\AspectUhatAl};
       \legend{{\Large Ti6Al4V}, {\Large SS316},{\Large AlSi10Mg}}
      \end{loglogaxis}
    \end{tikzpicture}
    }
  \caption{ Correlation of the aspect ratio with $\bm{Ma\hat{U}}= a_0\bm{Ma}+a_1\bm{MaE}+a_2\bm{MaPe}$, plotted on a log-log scale, for ~(\ref{fig:AspectPrMaUhatTi}) Ti6Al4V alloy, and for ~(\ref{fig:AspectPrMaUhatCombined}) three alloys (Ti6Al4V, SS316 and AlSi10Mg) shown in a single plot to demonstrate clustering.  \textcolor{black}{A combined plot demonstrating this clustering for all the five alloys (Ti6Al4V, SS316,  AlSi10Mg,  IN718 and AZ91D) considered in this work can be found in Figure~\ref{fig:AspectMaUhatCombinedSupplementary} of the Supplementary Information.} Each point in these plots represent a single simulation result for the relevant quantities plotted, and is obtained from the FEM framework.
}
  \label{fig:Aspect} 
\end{figure}
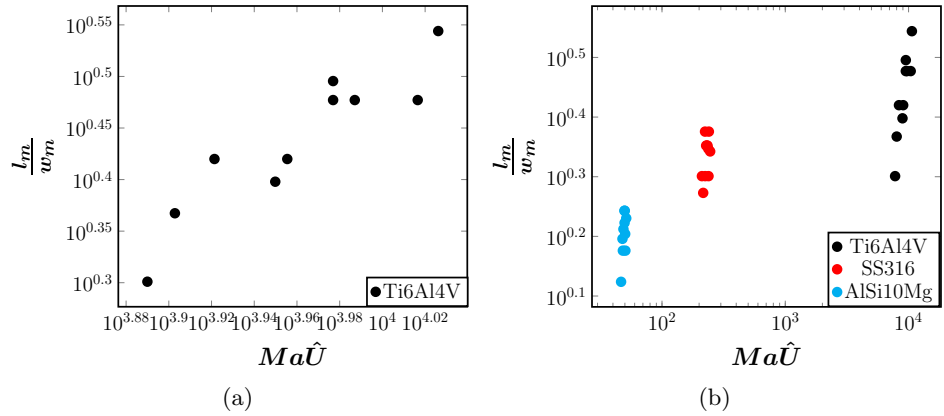

\subsection{Influence of Stefan number on meltpool volume}{\label{subsec:VolumeMeasure}}
\textcolor{black}{In this section, a relation between the input non-dimensional numbers and the volume of the meltpool is presented. }The volume of the meltpool is the volume of the localized region where the heat received by the material resulted in the \textcolor{black}{melting of the material. }\textcolor{black}{
the pool volume is identified using the liquid fraction, $\phi$, which is a variable that we track in the model at all time instances and all spatial points in the domain. }To completely melt the material, the material needs to absorb enough power to raise the temperature above the liquidus temperature and to overcome the latent heat barrier.  \textcolor{black}{
The dimensionless power absorbed by the material is associated as a form of latent power, $P_{L}=\frac{\Omega_m}{l_p^3}\bm{\frac{Tc}{Ste}}$. }Since $\bm{\frac{Tc}{Ste}}$ is constant for a given material,  $P_{L}$ is proportional to the volume of the material melted $(\Tilde{\Omega}_m)$.  One can expect, as a first-order argument, that more material will melt if $\hat{\bm{U}}$ is higher. Thus, one can expect $P_{L}=\Tilde{\Omega}_m\frac{\bm{Tc}}{\bm{Ste}} $ to increase with $ \hat{\bm{U}}$.  This understanding helps us anticipate that the meltpool volume for different alloys, $\tilde{\Omega}_m$, increases with $\frac{\bm{Ste\hat{U}}}{\bm{Tc}}$, and this can indeed be seen in Figure~\ref{fig:VolumeUhatAl}-\ref{fig:VolumeUhatTi}.  The numerical range of the meltpool volumes across the data points is similar, but for a given alloy, the meltpool volume increases with ${\frac{\bm{Ste}}{\bm{Tc}}\bm{\hat{U}}}$.  With this analysis,  an important correlation \textcolor{black}{is obtained} between the meltpool volume and the quantity $\bm{\frac{Ste}{Tc}\hat{U}} = a_0\frac{\bm{Ste}}{\bm{Tc}}+a_1\frac{\bm{SteE}}{\bm{Tc}}+a_2\frac{\bm{StePe}}{\bm{Tc}}$.

\pgfplotstableread{data/VolumeUhat.txt}{\VolUhat}
\pgfplotstableread{data/VolumeUhatSS.txt}{\VolUhatSS}
\pgfplotstableread{data/VolumeUhatTi.txt}{\VolUhatTi}
\pgfplotstableread{data/VolumeUhatAl.txt}{\VolUhatAl}

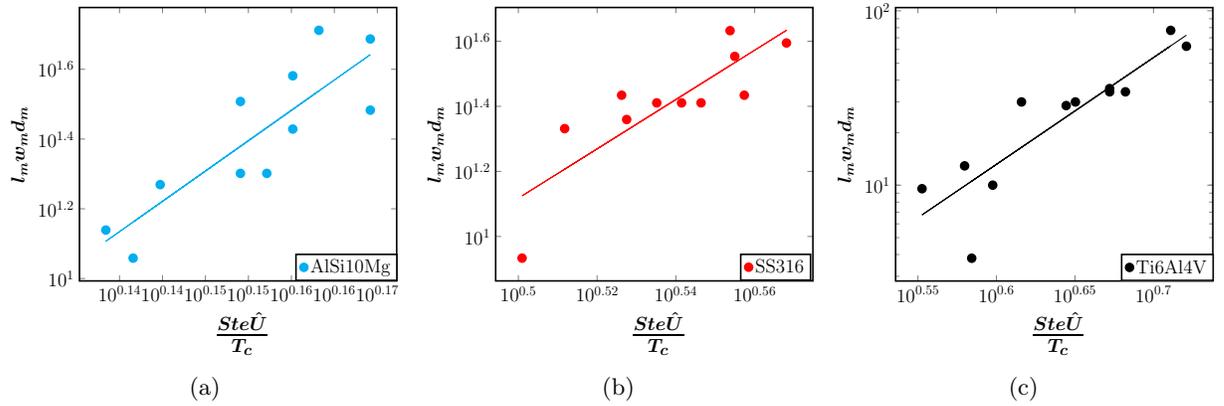
\begin{figure}[ht]
  \centering
  \subfloat[\label{fig:VolumeUhatAl}]{
    \begin{tikzpicture}[scale=0.5]
      \begin{loglogaxis}[ ticklabel style = {thick, font=\Large },very thick,xlabel={\huge {{ $\bm{\frac{Ste \hat{U}}{T_c}}$}}},ylabel={\Large \textbf{$\bm{l_m w_m d_m}$}},legend style={at={(0.68,0.1)}},x unit=, y unit=, mark repeat={1},scale only axis=true] 
    \addplot [only marks, cyan, mark=*, mark size=3] table [x={AL_x}, y={AL_Vol}] {\VolUhatAl};
        \addplot [thick, cyan] table[y={create col/linear regression={y=AL_Vol}}]{\VolUhatAl};  
        \legend{{\Large AlSi10Mg}}
      \end{loglogaxis}
    \end{tikzpicture} 
    }
    \subfloat[\label{fig:VolumeUhatSS}]{
    \begin{tikzpicture}[scale=0.5]
      \begin{loglogaxis}[ ticklabel style = {thick, font=\Large },very thick,xlabel={\huge {{ $\bm{\frac{Ste \hat{U}}{T_c}}$}}},ylabel={\Large \textbf{$\bm{l_m w_m d_m}$}},legend style={at={(0.77,0.1)}},x unit=, y unit=, mark repeat={1},scale only axis=true] 
       \addplot [only marks, red, mark=*, mark size=3] table [x={SS_x}, y={SS_Vol}] {\VolUhatSS};
         \addplot [thick, red] table[y={create col/linear regression={y=SS_Vol}}]{\VolUhatSS};  
        \legend{{\Large SS316}}
      \end{loglogaxis}
    \end{tikzpicture}
    } 
    \subfloat[\label{fig:VolumeUhatTi}]{
    \begin{tikzpicture}[scale=0.5]
      \begin{loglogaxis}[ ticklabel style = {thick, font=\Large },very thick,xlabel={\huge {{ $\bm{\frac{Ste \hat{U}}{T_c}}$}}},ylabel={\Large \textbf{$\bm{l_m w_m d_m}$}},legend style={at={(0.71,0.1)}},x unit=, y unit=, mark repeat={1},scale only axis=true] 
       \addplot [only marks, black, mark=*, mark size=3] table [x={TI_x}, y={TI_Vol}] {\VolUhatTi};
          \addplot [thick, black] table[y={create col/linear regression={y=TI_Vol}}]{\VolUhatTi}; 
        \legend{{\Large Ti6Al4V}}
      \end{loglogaxis}
    \end{tikzpicture}
    }
  \caption{Correlation of the meltpool volume $(l_mw_md_m)$ with $\bm{\frac{Ste \hat{U}}{T_c}}= a_0\frac{\bm{Ste}}{\bm{Tc}}+a_1\frac{\bm{SteE}}{\bm{Tc}}+a_2\frac{\bm{StePe}}{\bm{Tc}}$, plotted on a log-log scale, for (\ref{fig:VolumeUhatAl}) AlSi10Mg, (\ref{fig:VolumeUhatSS}) SS316, and (\ref{fig:VolumeUhatTi}) Ti6Al4V alloys.  \textcolor{black}{Corresponding plots comparing IN718 and AZ91D alloys, and a comparison of all the five alloys considered in this work can be found in Figure~\ref{fig:VolumeUhatSupplementary} and Figure~\ref{fig:VolumeUhatFull} of the Supplementary Information, respectively.} Each point in these plots represent a single simulation result for the relevant quantities plotted, and is obtained from the FEM framework.}
  \label{fig:VolumeUhat}
\end{figure}

\subsection{ \textcolor{black}{Influence of the heat absorbed on the solidification cooling rates}}{\label{subsec:MicroEvolu}}
\textcolor{black}{
In this section, a discussion on characterizing the microstructure of the solidified region is presented.  To support this discussion,  $\hat{\bm{U}}$ is used to explain the variation in the output variables like the temperature gradient, $G$, and the cooling rate, $G\bm{\nu}_p$, where $\bm{\nu}_p$ is the laser scan speed. }These variables are traditionally understood to have a direct influence on the microstructure in the solidification literature. $ \hat{\bm{U}}=a_0+a_1\bm{E}+a_2\bm{Pe}$, has a strong linear correlation with the temperature gradient, $G$, as shown in the Figures~\ref{fig:GUhatAl}-\ref{fig:GUhatTi}. Considering the correlations observed in these figures, \textcolor{black}{ it can be inferred} that the non-dimensional temperature gradient is proportional to the value of $\hat{\bm{U}}$. Further, this relation is expressed entirely in terms of input material properties and process parameters. 

\pgfplotstableread{data/GUhatSS.txt}{\GUhatSS}
\pgfplotstableread{data/GUhatTi.txt}{\GUhatTi}
\pgfplotstableread{data/GUhatAl.txt}{\GUhatAl}

\begin{figure}[ht] 
  \centering 
     \subfloat[b\label{fig:GUhatAl}]{ 
    \begin{tikzpicture}[scale=0.5]
      \begin{axis}[ ticklabel style = {thick, font=\Large },very thick,xlabel={\Large \textbf{$\hat{\bm{U}}=a_0+a_1\bm{E}+a_2\bm{Pe} $}},ylabel={\Large \textbf{$\tilde{G}_{max}$}},legend style={at={(0.68,0.1)}},x unit=, y unit=, mark repeat={1},scale only axis=true] 
    \addplot [only marks, cyan, mark=*, mark size=3] table [x={AL_Uhat}, y={AL_GM}] {\GUhatAl};
          \addplot [thick, cyan] table[y={create col/linear regression={y=AL_GM}}]{\GUhatAl}; 
       \legend{{\Large AlSi10Mg}}
      \end{axis}
    \end{tikzpicture}
    } 
       \subfloat[b\label{fig:GUhatSS}]{
    \begin{tikzpicture}[scale=0.5]
      \begin{axis}[ ticklabel style = {thick, font=\Large },very thick,xlabel={\Large \textbf{$\hat{\bm{U}}=a_0+a_1\bm{E}+a_2\bm{Pe} $}},ylabel={\Large \textbf{$\tilde{G}_{max}$}},legend style={at={(0.77,0.1)}},x unit=, y unit=, mark repeat={1},scale only axis=true] 
       \addplot [only marks, red, mark=*, mark size=3] table [x={SS_Uhat}, y={SS_GM}] {\GUhatSS};
         \addplot [thick, red] table[y={create col/linear regression={y=SS_GM}}]{\GUhatSS}; 
      \legend{{\Large SS316}}
      \end{axis}
    \end{tikzpicture}
    } 
       \subfloat[b\label{fig:GUhatTi}]{
    \begin{tikzpicture}[scale=0.5]
      \begin{axis}[ ticklabel style = {thick, font=\Large },very thick,xlabel={\Large \textbf{$\hat{\bm{U}}=a_0+a_1\bm{E}+a_2\bm{Pe} $}},ylabel={\Large \textbf{$\tilde{G}_{max}$}},legend style={at={(0.71,0.1)}},x unit=, y unit=, mark repeat={1},scale only axis=true] 
    \addplot [only marks, black, mark=*, mark size=3] table [x={TI_Uhat}, y={TI_GM}] {\GUhatTi};
       \addplot [thick, black] table[y={create col/linear regression={y=TI_GM}}]{\GUhatTi}; 
      \legend{{\Large Ti6Al4V}}
      \end{axis}
    \end{tikzpicture}
    }
  \caption{Dimensionless temperature gradient $(G)$ with the $\hat{\bm{U}}$ for different alloys.   \textcolor{black}{Plots corresponding to IN718 and AZ91D alloy material can be found in Figure~\ref{fig:DimlessGUhatSupplementary} of the Supplementary Information.} Each point in these plots represent a single simulation result for the relevant quantities plotted, and is obtained from the FEM framework.}
  \label{fig:DimlessGUhat}
\end{figure}
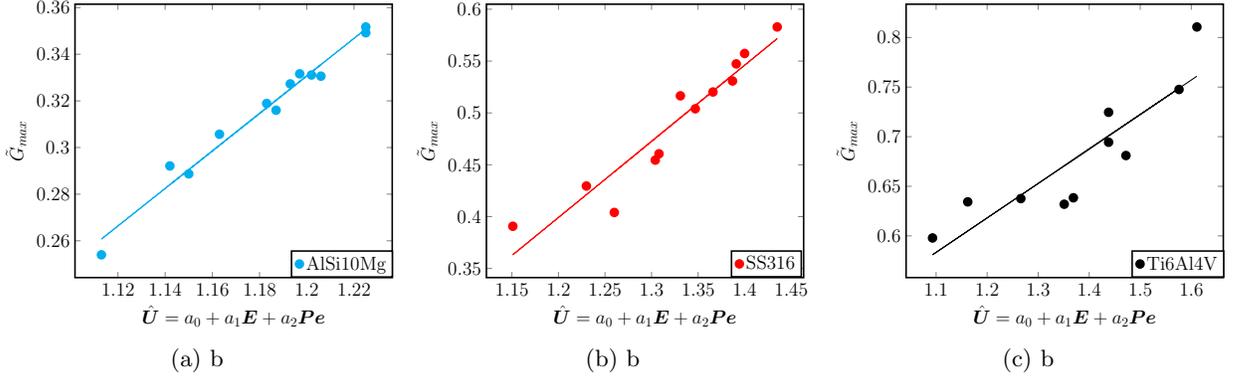


In the LPBF process, the laser scan speed controls the speed of movement of the solidification interface.  \textcolor{black}{The laser scan speed used in the simulations is assumed to be equal to the maximum solid-liquid interface velocity $(\bm{\nu}=\bm{\nu}_p $~\citep{mohammadpour2020revisiting, rappaz1990analysis}. }\textcolor{black}{The variation of the cooling rates, $G\bm{\nu}_p$ are observed with $\hat{\bm{U}}$ by increasing dimensionless laser power, $\bm{E}$, but keeping the P\'eclet number fixed.  The dimensionless temperature gradient, $\tilde{G}_{max}$, will increase with $\hat{\bm{U}}$. 
as seen in the Figure~\ref{fig:DimlessGUhat}}. It is also instructional to see that the cooling rate, $G\bm{{\nu}}_p$, increases with $\hat{\bm{U}}$, as shown in Figure~\ref{subfig:GVwithUhat1}. There exists a well-known correlation between a microstructure size, $\lambda_2$, and the cooling rate, $G\bm{\nu}_p$, in the solidification process, and is given by the relation $\lambda_2=25(G\bm{\nu}_p)^{-0.28}$~\citep{thoma1995directed}.  Using this relation, \textcolor{black}{it can be seen that} the size of the microstructure will get finer as we increase $\hat{\bm{U}}$, which can be achieved by increasing $\bm{E}$. Thus, the size of the microstructure correlates with the input non-dimensional quantity given by $\hat{\bm{U}}=a_0+a_1\bm{E}+a_2\bm{Pe}$. 

\textcolor{black}{
The effect of $\hat{\bm{U}}$ on the cooling rate is studied due to a change in the P\'eclet number, $\bm{Pe}$, but keeping the dimensionless laser power, $\bm{E}$, fixed.  $\hat{\bm{U}}$ decreases with increase in $\bm{Pe}$ and $\bm{E}$ fixed.  From Figure~\ref{fig:DimlessGUhat},  it is known that the dimensionless temperature gradient, $G$, and its dimensional counterpart, both decrease if we decrease $\hat{\bm{U}}$. }Thus the cooling rate, $G\bm{\nu}_p$, decreases with increase in  $\hat{\bm{U}}$, as seen in Figure~\ref{subfig:GVwithUhat2}. This decrease is solely due to an increase in the value of $\bm{Pe}$. Using this information,  a correlation \textcolor{black}{is identified} between the change in the size of the microstructure and the P\'eclet number, i.e, increasingly coarser microstructural features can be observed with a decrease in the cooling rate.

\pgfplotstableread{data/Microstruct1.txt}{\MSone}
\pgfplotstableread{data/Microstruct2.txt}{\MStwo}
\pgfplotstableread{data/Microstruct3.txt}{\MSthree}

\begin{figure}[ht]
  \centering 
\subfloat[\label{subfig:GVwithUhat1}]{
   \begin{tikzpicture}[scale=0.55]
      \begin{axis}[ ticklabel style = {thick, font=\Large },very thick,xlabel={\LARGE { $\hat{\bm{U}}$ }},ylabel={\LARGE \textbf{$G\bm{\nu}_p$}},legend style={at={(0.4,0.33)}},x unit=, y unit=, mark repeat={1},scale only axis=true] 
       \addplot [only marks, red, mark=*, mark size=3] table [y={GV_3}, x={Uhat}] {\MSthree};
      \end{axis}
    \end{tikzpicture}
  }
\subfloat[\label{subfig:GVwithUhat2}]{
   \begin{tikzpicture}[scale=0.55]
    \begin{axis}[ ticklabel style = {thick, font=\Large },very thick,xlabel={\LARGE { $\hat{\bm{U}}$ }},ylabel={\LARGE \textbf{$G\bm{\nu}_p$}},legend style={at={(0.4,0.33)}},x unit=, y unit=, mark repeat={1},scale only axis=true]
       \addplot [only marks, red, mark=*, mark size=3] table [y={GV_2}, x={Uhat}] {\MStwo};
      \end{axis}
    \end{tikzpicture}
  }

  \caption{~(\ref{subfig:GVwithUhat1}) Variation of the dimensional cooling rate, $G\bm{\nu}_p$, with $\bm{\hat{U}}= a_0+a_a\bm{E}+a_2\bm{Pe}$, plotted on a log-log scale. Here, $\bm{\hat{U}}$ is changed by changing $\bm{E}$, but keeping $\bm{Pe}$ fixed for SS316 alloy,  ~(\ref{subfig:GVwithUhat2}) Variation of dimensional cooling rate, $G\bm{\nu}_p$, with $\hat{\bm{U}}=a_0+a_a\bm{E}+a_2\bm{Pe}$, plotted on a log-log scale. Here, $\bm{\hat{U}}$ is changed by changing $\bm{Pe}$, but keeping $\bm{E}$ fixed for SS316 alloy.}
  \label{fig:Microstructtype1}
\end{figure}
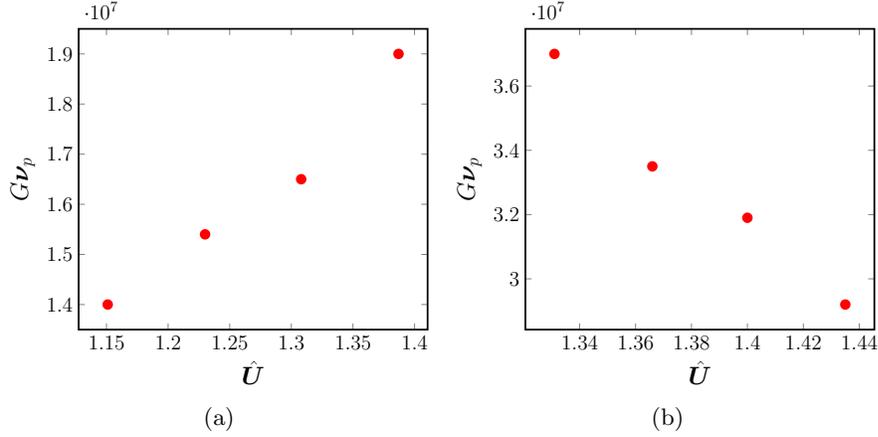

\section{Conclusion}{\label{sec:Conclusions}}
\textcolor{black}{ 
In this work, the meltpool dynamics of the laser powder bed fusion (LPBF) process are numerically modeled and connections are made to important dimensionless quantities influencing the thermo-fluidic evolution of the meltpool and its morphology. Processes like the interaction of the moving laser power source with the powdered metal, formation of the meltpool, its subsequent solidification, etc., make LPBF a highly coupled multiphysics process. To investigate the multiphysics interactions, the thermo-fluidic governing equations relevant to the LPBF process are numerically modeled using a Finite Element Method (FEM) framework. The simulation predictions were validated by comparing with available results from the literature and with experimental observations of the cooling rates available from our experimental collaborators. Using the classical Buckingham-$\pi$ theorem and a careful choice of relevant characteristic quantities, the governing equations were reduced to their dimensionless form. Using the dimensionless form and the FEM simulations, an important dimensionless quantity, interpreted as the heat absorbed by the metal powder and the meltpool, is identified. Around sixty different cases of the LPBF process were simulated by varying the alloy type and the process conditions, and the simulation data was used to obtain an explicit form of the dimensionless heat absorbed in terms of the input dimensionless numbers using the method of linear least-squares regression. Using physics-based and statistical arguments, a linear model showing dependence of the heat absorbed on the P\'eclet number and the dimensionless power is established. }

\textcolor{black}{ 
The measure of advection inside the meltpool is quantified in terms of the P\'eclet and the Marangoni numbers and it is found that materials such as Ti-6Al-4V show greater advection, represented by $\bm{Pe^{*}}$, and an elongated elliptical meltpool. Materials like AlSi10Mg show the least advection and whereas SS316 shows moderate amount of advection. It was found that the meltpool volume of materials such as Ti-6Al-4V, SS316, and AlSi10Mg increases with the product of Stefan number and dimensionless measure of heat absorbed. Solidification cooling rates decrease with the measure of heat absorbed if the P\'eclet number is reduced keeping dimensionless power fixed. This characterization of the meltpool morphology using classical dimensionless numbers and the impact of dimensionless power and P\'eclet number on the solidification cooling rates is a novel contribution of this work.}

\textcolor{black}{
In a future work, potential extension of this dimensional analysis framework to investigate meltpool characteristics such as keyhole formation and microstructural features such as grain morphology in the solidified meltpool region will be explored.}
 
\section*{Acknowledgement}
The authors would like to thank Prof. Dan Thoma (University of Wisconsin-Madison) and Dr. Kaila Bertsch (University of Wisconsin-Madison; now at Lawrence Livermore National Laboratory) for very useful discussions on microstructure evolution during the LPBF process, and for providing the experimental data on cooling rates that are shown in Figure~\ref{fig:Experiment_Numerical}.

\bibliography{main}

\clearpage
\setcounter{table}{0} 
\setcounter{figure}{0} 
\setcounter{section}{0} 
\setcounter{page}{1} 

\section*{Supplementary Information for {\normalfont ``A numerical investigation of dimensionless numbers characterizing meltpool morphology of the laser powder bed fusion process''}}{\label{sec:SupplementaryInformation}}
Kunal Bhagat, Shiva Rudraraju\\
Department of Mechanical Engineering, University of Wisconsin-Madison, Madison, WI, USA\\

\vspace{0.1in}

\subsection*{Material properties and process variables}
Material properties of various alloys and the corresponding AM process variables used in this work are listed here. These properties were collected from multiple sources in the literature.  The properties of solid and liquid materials are averaged and dependence on the temperature is neglected. 
\begin{table}[ht]
\centering 
\caption{Average material properties for different alloys used to calculate input non-dimensional numbers of the thermo-fluidic model ~\cite{mukherjee2018heat}, \cite{shen2020thermo}, \cite{Properties}} 
\begin{tabular}{c c c c c c c} 
\hline 
Property & SS316 & Ti6Al4V & IN718 & AlSi10Mg & AZ91D  \\ [0.5ex] 
\hline 
$\rho (\frac{kg}{m^3} )$& 7800  & 4000& 8100 &2670 & 1675   \\[0.5ex] 
$c(\frac{J}{kgK}) $ & 490  & 570 & 435 & 890 & 1122  \\[0.5ex]
$k(\frac{W}{mk})$ & 36.5 & 7.3 & 11.4  & 173.0 & 77 .5  \\[0.5ex]
$\mu (\frac{Kg}{ms})$ & $7.0\times10^{-3}$ & $4.0\times10^{-3} $& $5.0\times10^{-3} $& $1.3\times10^{-3}$& $3.0\times10^{-3}$   \\[0.5ex]
$\frac{d\gamma}{dT}(\frac{N}{mK})$    & $-4.00\times10^{-4}$ & $-2.63\times10^{-3} $& $-3.70\times10^{-3} $& $-3.5\times10^{-4}$& $-2.13\times10^{-4}$ \\[0.5ex]
$\beta (\frac{1}{K})$ & $5.85\times10^{-5}$ & $2.50\times10^{-5} $& $4.8\times10^{-5} $& $2.4\times10^{-5}$& $9.54\times10^{-5}$    \\[0.5ex]
$\kappa (m^2)$ & $5.56\times10^{-13}$ & $5.56\times10^{-13} $& $5.56\times10^{-13}$& $5.56\times10^{-13}$& $5.56\times10^{-13}$    \\[0.5ex]
$L(\frac{J}{kg})$ & $2.72\times10^{5}$ & $2.84\times10^{5} $& $2.09\times10^{5}$& $4.23\times10^{5}$& $3.73\times10^{5}$    \\[0.5ex]
$T_{s} ($K$ )$ & 1693 &1878 &1533 & 831&  743 \\[0.5ex]
$T_{l} ($K$ )$ & 1733 & 1928 &1609 & 867&  868\\[0.5ex]
\hline 
\end{tabular}
\label{table:AllMaterialProperties} 
\end{table}

\begin{table}[ht]  
\centering 
\caption{Chosen process conditions for different alloys used to calculate input non-dimensional numbers of the thermo-fluidic model} 
\begin{tabular}{c c c c c c c} 
\hline 
Material &  (Laser Power,  Scan speed)  $ (P,\bm{\nu}_p)$   \\ [0.5ex] 
\hline 
&(70, 0.3),  (80, 0.4),  (90, 0.5),  (100, 0.6) \\ 
SS316&\quad  (110, 0.7),  (110, 0.8),  (110, 0.9),  (110, 1.0) \\ 
&(65, 0.5),  (75, 0.5),  (85, 0.5),  (95, 0.5)\\[1.0ex] 
 \hline 
&(15, 0.2),  (25, 0.5),  (35, 0.7),  (45, 0.9)\\ 
Ti6Al4V&(40, 0.6),  (40, 0.7),  (40, 0.8),  (40, 1.0)\\ 
&(35, 0.9),  (40, 0.9),  (45, 0.9),  (50, 0.9)\\[1.0ex] 
\hline 
&\quad (20, 0.15),  (30, 0.25),  (40, 0.45),  (50, 0.75),\\ 
IN718&(45, 0.80),  (45, 0.90),  (45, 1.0),  (45, 1.1)\\ 
&\quad (53, 0.95),  (55, 0.95),  (58, 0.95),  (60, 0.95)\\[1.0ex] 
\hline 
&\quad(75, 0.35),  (85, 0.45),  (95, 0.55),  (105, 0.65) \\ 
AlSi10Mg&\quad (100, 0.6),  (100, 0.7),  (100, 0.8),  (100, 0.9)\\ 
&(90, 1.1),  (95, 1.1),  (100, 1.1),  (110, 1.1)\\[1.0ex] 
\hline 
&(35, 0.25),  (40, 0.30),  (45, 0.35),  (50, 0.45\\ 
AZ91D&(40, 0.30),  (40, 0.40),  (40, 0.50),  (40, 0.60)\\ 
&(40, 0.60),  (50, 0.60),  (60, 0.60),  (70, 0.60)\\ 
\hline 
\end{tabular}
\label{table:AllProcessProperties} 
\end{table}

\clearpage 
\subsection*{Additional correlations of the dimensionless numbers}

\textcolor{black}{
The influence of the P\'eclet number on advection transport in the meltpool for additional materials IN718 and AZ91D is shown in Figure~\ref{fig:AdvectMeasureSupplementary}, and the corresponding comparison of the five alloys considered in this work is given in Figure~\ref{fig:AdvectMeasureFull}.  The combined plot of aspect ratio with $\bm{Ma}\bm{\hat{U}}$ is given in the \textcolor{black}{Figure}~\ref{fig:AspectMaUhatCombinedSupplementary}. The meltpool aspect ratio and Marangoni number of IN718 and Ti6AL4V are similar magnitudes.  Similarly, the meltpool aspect ratio and Marangoni number of AZ91D and AlSi10Mg are comparable.  The influence of Stefan number on the meltpool volumes for IN718 and AZ91D is shown in Figure~\ref{fig:VolumeUhatSupplementary}, and the corresponding comparison of the five alloys considered in this work is given in Figure~\ref{fig:VolumeUhatFull}. Influence of dimensionless heat absorbed on  \textcolor{black}{non-dimensional temperature gradient ($G$)} for IN718 and AZ91D is shown in Figure~\ref{fig:DimlessGUhatSupplementary}.}

\pgfplotstableread{data/AdvectMeasureIN.txt}{\AdvectMIN}
\pgfplotstableread{data/AdvectMeasureAZ.txt}{\AdvectMAZ}

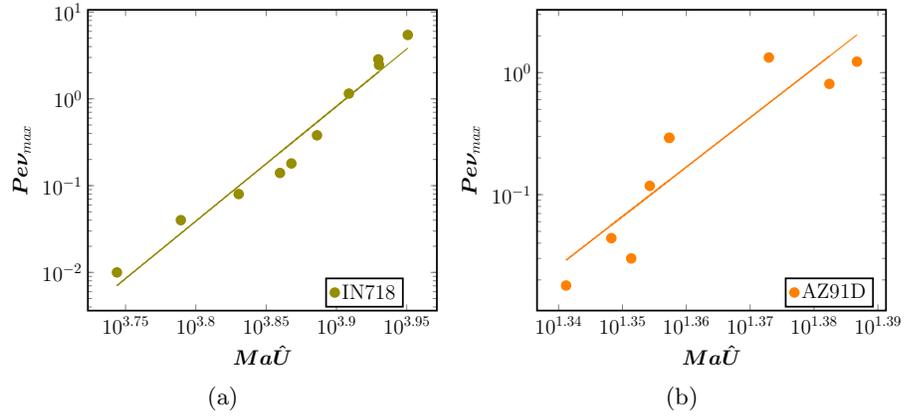
\begin{figure}[h]
  \centering
\subfloat[\label{AdvectIN}]{
   \begin{tikzpicture}[scale=0.55]
      \begin{loglogaxis}[ ticklabel style = {thick, font=\Large },very thick,xlabel={\Large { $\bm{Ma\hat{U}}$ }},ylabel={\Large \textbf{$\bm{Pe}\bm{\nu}_{max}$}},legend style={at={(0.68,0.11)}},x unit=, y unit=, mark repeat={1},scale only axis=true] 
     \addplot [only marks, olive, mark=*, mark size=3] table [x={IN_MaUhat}, y={IN_PeU0}] {\AdvectMIN};   
   \addplot [thick, olive] table[y={create col/linear regression={y=IN_PeU0}}]{\AdvectMIN};  
      \legend{{\Large IN718}}
      \end{loglogaxis}
    \end{tikzpicture}
  }
    \subfloat[\label{AdvectAZ}]{
   \begin{tikzpicture}[scale=0.55]
      \begin{loglogaxis}[ ticklabel style = {thick, font=\Large },very thick,xlabel={\Large { $\bm{Ma\hat{U}}$ }},ylabel={\Large \textbf{$\bm{Pe}\bm{\nu}_{max}$}},legend style={at={(0.71,0.11)}},x unit=, y unit=, mark repeat={1},scale only axis=true] 
       \addplot [only marks, orange, mark=*, mark size=3] table [x={AZ_MaUhat}, y={AZ_PeU0}] {\AdvectMAZ};
          \addplot [thick, orange] table[y={create col/linear regression={y=AZ_PeU0}}]{\AdvectMAZ};  
             \legend{{\Large AZ91D}}
      \end{loglogaxis}
    \end{tikzpicture}
  }
  \caption{ \textcolor{black}{Measure of total advection measured as $\mathbf{Pe\bm{\nu}_{max}}$ vs surface tension based advection $\mathbf{Ma\hat{U}}=a_0\bm{Ma} +a_1\bm{MaE}+a_2\bm{MaPe}$ on a log-log scale for (\ref{AdvectIN})IN718  (\ref{AdvectAZ})AZ91D, alloys.}}
  \label{fig:AdvectMeasureSupplementary} 
\end{figure}

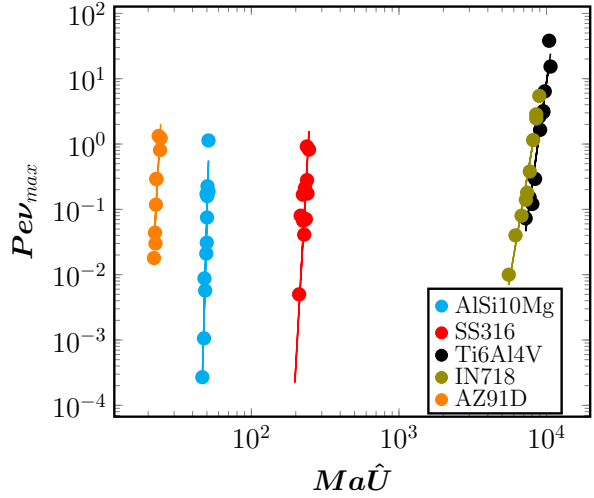
\begin{figure}[ht]
  \centering
   \begin{tikzpicture}[scale=0.75]
      \begin{loglogaxis}[ ticklabel style = {thick, font=\Large },very thick,xlabel={\Large { $\bm{Ma\hat{U}}$ }},ylabel={\Large \textbf{$\bm{Pe}\bm{\nu}_{max}$}},legend style={at={(0.8,0.11)}},x unit=, y unit=, mark repeat={1},scale only axis=true] 
     \addplot [only marks, cyan, mark=*, mark size=3] table [x={AL_MaUhat}, y={AL_PeU0}] {\AdvectMAl}; \label{p11}
   \addplot [thick, cyan] table[y={create col/linear regression={y=AL_PeU0}}]{\AdvectMAl};  
       \addplot [only marks, red, mark=*, mark size=3] table [x={SS_MaUhat}, y={SS_PeU0}] {\AdvectMSS}; \label{p22}
        \addplot [thick, red] table [x={Logx}, y={FitLogy}] {\AdvectMSS};
        \addplot [only marks, black, mark=*, mark size=3] table [x={TI_MaUhat}, y={TI_PeU0}] {\AdvectMTi}; \label{p33}
        \addplot [thick, black] table[y={create col/linear regression={y=TI_PeU0}}]{\AdvectMTi};  
        \addplot [only marks, olive, mark=*, mark size=3] table [x={IN_MaUhat}, y={IN_PeU0}] {\AdvectMIN};   \label{p44}
   	\addplot [thick, olive] table[y={create col/linear regression={y=IN_PeU0}}]{\AdvectMIN};  
	\addplot [only marks, orange, mark=*, mark size=3] table [x={AZ_MaUhat}, y={AZ_PeU0}] {\AdvectMAZ}; \label{p55}
        \addplot [thick, orange] table[y={create col/linear regression={y=AZ_PeU0}}]{\AdvectMAZ};  
\node [draw,fill=white] at (rel axis cs: 0.8,0.16) {\shortstack[l]{
\ref{p11} \large AlSi10Mg \\
\ref{p22} \large SS316 \\
\ref{p33} \large Ti6Al4V \\ 
\ref{p44} \large IN718 \\
\ref{p55} \large AZ91D }};
\end{loglogaxis}
\end{tikzpicture}
  \caption{ Measure of total advection measured as $\mathbf{Pe\bm{\nu}_{max}}$ vs surface tension based advection $\mathbf{Ma\hat{U}}=a_0\bm{Ma} +a_1\bm{MaE}+a_2\bm{MaPe}$, plotted on a log-log scale, for all the five alloys (AlSi10Mg,  SS316, Ti6Al4V,  IN718 and AZ91D) considered in this work. }
  \label{fig:AdvectMeasureFull} 
\end{figure}

\pgfplotstableread{data/AspectRatioSS.txt}{\AspectUhatSS}
\pgfplotstableread{data/AspectRatioTi.txt}{\AspectUhatTi}
\pgfplotstableread{data/AspectRatioAl.txt}{\AspectUhatAl}
\pgfplotstableread{data/AspectRatioIN.txt}{\AspectUhatIN}
\pgfplotstableread{data/AspectRatioAZ.txt}{\AspectUhatAZ}

\begin{figure}[ht]
  \centering    
    \begin{tikzpicture}[scale=0.75]
      \begin{loglogaxis}[ ticklabel style = {thick, font=\Large },very thick,xlabel={\LARGE { $\bm{{Ma\hat{U}}}$ }},ylabel={\huge\bm{$\frac{l_{m}}{w_{m}}$}},legend style={at={(0.48,0.36)}},x unit=, y unit=, mark repeat={1},scale only axis=true] 
         \addplot [only marks, black, mark=*, mark size=3] table [x={TI_MaU}, y={TI_AR}] {\AspectUhatTi};
       \addplot [only marks, red, mark=*, mark size=3] table [x={SS_MaU}, y={SS_AR}] {\AspectUhatSS}; 
    \addplot [only marks, cyan, mark=*, mark size=3] table [x={AL_MaU}, y={AL_AR}] {\AspectUhatAl};
     \addplot [only marks, olive, mark=*, mark size=3] table [x={IN_MaU}, y={IN_AR}] {\AspectUhatIN};
    \addplot [only marks, orange, mark=*, mark size=3] table [x={AZ_MaU}, y={AZ_AR}] {\AspectUhatAZ};
       \legend{{\large Ti6Al4V}, {\large SS316},{\large AlSi10Mg}, {\large IN718},{\large AZ91D}}
      \end{loglogaxis}
    \end{tikzpicture}
  \caption{ \textcolor{black}{Correlation of the aspect ratio with $\bm{Ma\hat{U}}= a_0\bm{Ma}+a_1\bm{MaE}+a_2\bm{MaPe}$, plotted on a log-log scale, for all five alloys (Ti6Al4V, SS316,  AlSi10Mg,  IN718 and AZ91D) shown in a single plot to demonstrate clustering. }}
  \label{fig:AspectMaUhatCombinedSupplementary} 
\end{figure}
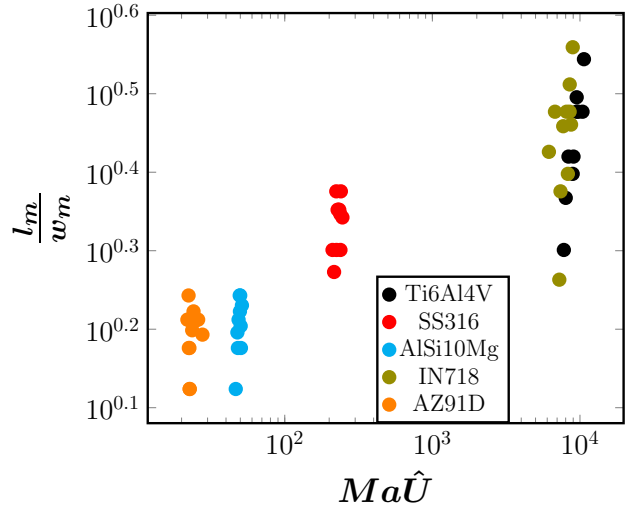

\pgfplotstableread{data/VolumeUhatIN.txt}{\VolUhatIN}
\pgfplotstableread{data/VolumeUhatAZ.txt}{\VolUhatAZ}
\begin{figure}[ht]
  \centering
  \subfloat[\label{fig:VolumeUhatIN}]{
    \begin{tikzpicture}[scale=0.5]
      \begin{loglogaxis}[ ticklabel style = {thick, font=\Large },very thick,xlabel={\huge {{ $\bm{\frac{Ste \hat{U}}{T_c}}$}}},ylabel={\Large \textbf{$\bm{l_m w_m d_m}$}},legend style={at={(0.68,0.1)}},x unit=, y unit=, mark repeat={1},scale only axis=true] 
    \addplot [only marks, olive, mark=*, mark size=3] table [x={IN_x}, y={IN_Vol}] {\VolUhatIN};
        \addplot [thick, olive] table[y={create col/linear regression={y=IN_Vol}}]{\VolUhatIN};  
        \legend{{\Large IN718}}
      \end{loglogaxis}
    \end{tikzpicture}
    }
    \subfloat[\label{fig:VolumeUhatAZ}]{
    \begin{tikzpicture}[scale=0.5]
      \begin{loglogaxis}[ ticklabel style = {thick, font=\Large },very thick,xlabel={\huge {{ $\bm{\frac{Ste \hat{U}}{T_c}}$}}},ylabel={\Large \textbf{$\bm{l_m w_m d_m}$}},legend style={at={(0.71,0.1)}},x unit=, y unit=, mark repeat={1},scale only axis=true] 
       \addplot [only marks, orange, mark=*, mark size=3] table [x={AZ_x}, y={AZ_Vol}] {\VolUhatAZ};
          \addplot [thick, orange] table[y={create col/linear regression={y=AZ_Vol}}]{\VolUhatAZ}; 
        \legend{{\Large AZ91D}}
      \end{loglogaxis}
    \end{tikzpicture}
    }
  \caption{\textcolor{black}{Correlation of the meltpool volume $(l_mw_md_m)$ with $\bm{\frac{Ste \hat{U}}{T_c}}= a_0\frac{\bm{Ste}}{\bm{Tc}}+a_1\frac{\bm{SteE}}{\bm{Tc}}+a_2\frac{\bm{StePe}}{\bm{Tc}}$, plotted on a log-log scale, for (\ref{fig:VolumeUhatIN}) IN718, and (\ref{fig:VolumeUhatAZ}) AZ91D, alloys.}}
  \label{fig:VolumeUhatSupplementary}
\end{figure}
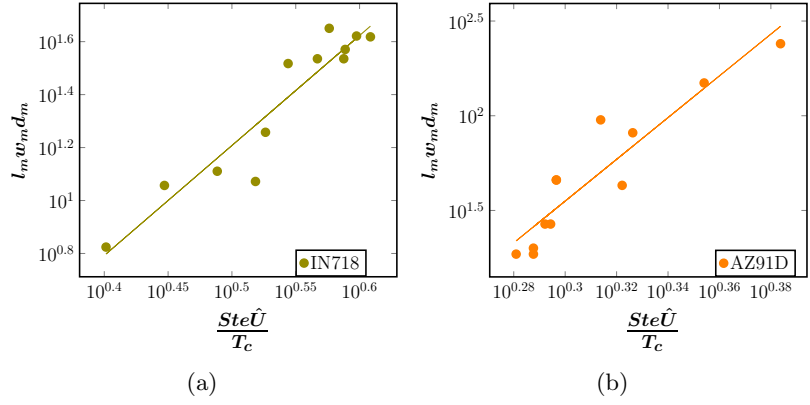

\begin{figure}[ht]
  \centering
    \begin{tikzpicture}[scale=0.85]
      \begin{loglogaxis}[ ticklabel style = {thick, font=\Large },very thick,xlabel={\huge {{ $\bm{\frac{Ste \hat{U}}{T_c}}$}}},ylabel={\Large \textbf{$\bm{l_m w_m d_m}$}},legend style={at={(0.68,0.1)}},x unit=, y unit=, mark repeat={1},scale only axis=true] 
   \addplot [only marks, cyan, mark=*, mark size=3] table [x={AL_x}, y={AL_Vol}] {\VolUhatAl}; \label{p111}
   \addplot [thick, cyan] table[y={create col/linear regression={y=AL_Vol}}]{\VolUhatAl};  
   \addplot [only marks, red, mark=*, mark size=3] table [x={SS_x}, y={SS_Vol}] {\VolUhatSS}; \label{p222}
   \addplot [thick, red] table[y={create col/linear regression={y=SS_Vol}}]{\VolUhatSS};  
   \addplot [only marks, black, mark=*, mark size=3] table [x={TI_x}, y={TI_Vol}] {\VolUhatTi}; \label{p333}
   \addplot [thick, black] table[y={create col/linear regression={y=TI_Vol}}]{\VolUhatTi};  
    \addplot [only marks, olive, mark=*, mark size=3] table [x={IN_x}, y={IN_Vol}] {\VolUhatIN}; \label{p444}
    \addplot [thick, olive] table[y={create col/linear regression={y=IN_Vol}}]{\VolUhatIN};  
    \addplot [only marks, orange, mark=*, mark size=3] table [x={AZ_x}, y={AZ_Vol}] {\VolUhatAZ}; \label{p555}
    \addplot [thick, orange] table[y={create col/linear regression={y=AZ_Vol}}]{\VolUhatAZ}; 
\node [draw,fill=white] at (rel axis cs: 0.3,0.16) {\shortstack[l]{
\ref{p111} \large AlSi10Mg \\
\ref{p222} \large SS316 \\
\ref{p333} \large Ti6Al4V \\ 
\ref{p444} \large IN718 \\
\ref{p555} \large AZ91D }};
\end{loglogaxis}
\end{tikzpicture}
  \caption{Correlation of the meltpool volume $(l_mw_md_m)$ with $\bm{\frac{Ste \hat{U}}{T_c}}= a_0\frac{\bm{Ste}}{\bm{Tc}}+a_1\frac{\bm{SteE}}{\bm{Tc}}+a_2\frac{\bm{StePe}}{\bm{Tc}}$, plotted on a log-log scale, for all the five alloys (AlSi10Mg,  SS316, Ti6Al4V,  IN718 and AZ91D) considered in this work. }
  \label{fig:VolumeUhatFull} 
\end{figure}
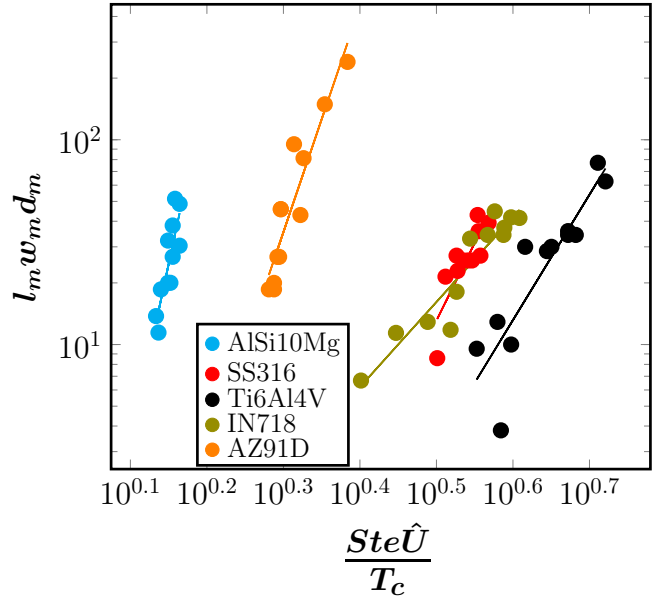

\pgfplotstableread{data/GUhatIN.txt}{\GUhatIN}
\pgfplotstableread{data/GUhatAZ.txt}{\GUhatAZ}

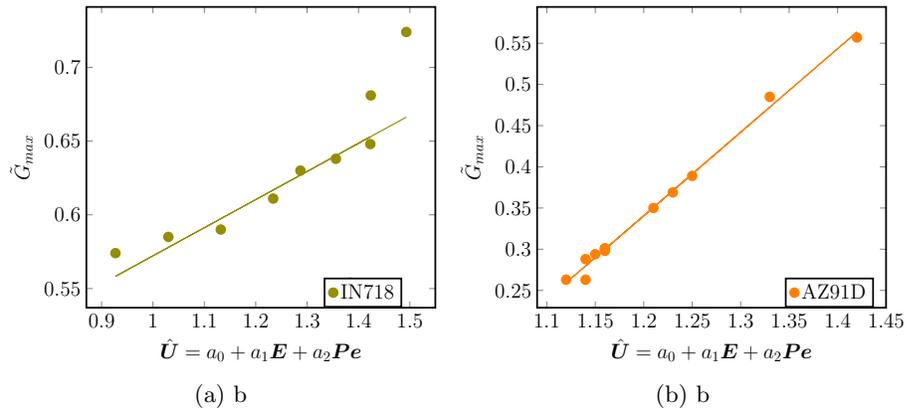
\begin{figure}[ht]
  \centering
     \subfloat[b\label{fig:GUhatIN}]{
    \begin{tikzpicture}[scale=0.55]
      \begin{axis}[ ticklabel style = {thick, font=\Large },very thick,xlabel={\Large \textbf{$\hat{\bm{U}}=a_0+a_1\bm{E}+a_2\bm{Pe} $}},ylabel={\Large \textbf{$\tilde{G}_{max}$}},legend style={at={(0.68,0.1)}},x unit=, y unit=, mark repeat={1},scale only axis=true] 
    \addplot [only marks, olive, mark=*, mark size=3] table [x={IN_Uhat}, y={IN_GM}] {\GUhatIN};
          \addplot [thick, olive] table[y={create col/linear regression={y=IN_GM}}]{\GUhatIN}; 
       \legend{{\Large IN718}}
      \end{axis}
    \end{tikzpicture}
    }
       \subfloat[b\label{fig:GUhatAZ}]{
    \begin{tikzpicture}[scale=0.55]
      \begin{axis}[ ticklabel style = {thick, font=\Large },very thick,xlabel={\Large \textbf{$\hat{\bm{U}}=a_0+a_1\bm{E}+a_2\bm{Pe} $}},ylabel={\Large \textbf{$\tilde{G}_{max}$}},legend style={at={(0.71,0.1)}},x unit=, y unit=, mark repeat={1},scale only axis=true] 
    \addplot [only marks, orange, mark=*, mark size=3] table [x={AZ_Uhat}, y={AZ_GM}] {\GUhatAZ};
       \addplot [thick, orange] table[y={create col/linear regression={y=AZ_GM}}]{\GUhatAZ}; 
      \legend{{\Large AZ91D}}
      \end{axis}
    \end{tikzpicture}
    }
  \caption{\textcolor{black}{Dimensionless temperature gradient $(G)$ with the $\hat{\bm{U}}$ for (\ref{fig:GUhatIN}) IN718 and (\ref{fig:GUhatAZ}) AZ91D alloys.}}
  \label{fig:DimlessGUhatSupplementary}
\end{figure}

\end{document}